\documentclass[leqno,11pt]{article}

\usepackage{amssymb}
\usepackage{amsfonts}
\usepackage{amscd}
\usepackage{amstext}
\usepackage{amssymb}
\usepackage{latexsym}
\usepackage{amsmath}
\usepackage{palatino}
\usepackage{color}

\newtheorem{theorem}{Theorem}[section]

\newtheorem{lemma}[theorem]{Lemma}

\newtheorem{claim}[theorem]{Claim}

\newcommand{\R}{\mathbb R}
\newcommand{\N}{\mathbb N}

\date{}

\begin{document}


\begin{center}
		{\Large Existence and concentration of nontrivial solitary waves for a generalized Kadomtsev--Petviashvili equation in $\mathbb{R}^2$}
\end{center}
	
\begin{center}
		{\small
			Claudianor O. Alves \\ Unidade Acad\^emica de Matem\'atica \\ Universidade Federal de Campina Grande \\ 58429-970, Campina Grande - PB, Brazil
			\\ E-mail: {\tt coalves@mat.ufcg.edu.br}
		}
\end{center}
	
\vspace{-0,3cm}
	
\begin{center}
		{\small
			Chao Ji$^{*}$ \\ School of Mathematics \\
			East China University of Science and Technology\\ Shanghai, 200237,  China \\
			E-mail: {\tt jichao@ecust.edu.cn}
		}
\end{center}
	
\vspace{-0,2cm}

\noindent {\sc Abstract}.  In this paper, we study the existence and concentration of solitary waves for
a class of generalized Kadomtsev-Petviashvili equations with the potential in $\mathbb{R}^2$ via the variational methods.

	
\noindent{\sc AMS Subject Classification MSC 2020}. 35A15, 35A18, 76B25.

\vspace{-0.5cm}
	
\footnotetext{$^{*}$ Corresponding author.}

\section{Introduction}

This paper is devoted to existence of the nontrivial solitary waves for
the generalized Kadomtsev-Petviashvili equation with a potential $\widetilde{V}$ and a parameter $\varepsilon>0$, namely
$$
v_t +\varepsilon^2 v_{xxx} + \left(f(v)- \widetilde{V}(x,y)v\right)_{x}=D^{-1}_{x}v_{yy},
\,\, \mbox{in} \,\,  \mathbb{R}^{+}\times\mathbb{R}^2\,\,  \eqno{(\widetilde{P}_{\varepsilon})}
$$
where
$$
D_{x}^{-1} h(x, y)=\int_{-\infty}^{x} h(s, y)ds.
$$

A solitary wave is a solution of the form
$v(t, x, y) = u(x-ct, y)$, $c>0$ is fixed. Substituting in $(\widetilde{P}_{\varepsilon})$, we obtain
$$
\varepsilon^2 u_{xxx}-cu_{x}+ \left(f(u)- \widetilde{V}(x,y)u\right)_{x}=D^{-1}_{x}u_{yy},
\,\, \mbox{in} \,\,  \mathbb{R}^2,
$$
or equivalently,
$$
\Big(- \varepsilon^2 u_{xx} +D^{-2}_{x}u_{yy}+   V(x,y)u-f(u)\Big)_{x}=0,
\,\, \mbox{in} \,\,  \mathbb{R}^2, \eqno{(P_{\varepsilon}^{*})}
$$
where $V=\widetilde{V} +c$, and after a scaling one obtains
$$
\left(-  u_{xx} +D^{-2}_{x}u_{yy}+   V(\varepsilon x,\varepsilon y)u-f(u)\right)_{x}=0,
\,\, \mbox{in} \,\,  \mathbb{R}^2. \eqno{(P_{\varepsilon})}
$$
If $V(\varepsilon x,\varepsilon y)=\tau $ is a positive constant, problem $(P_{\varepsilon})$ in $\mathbb{R}^N$ becomes
\begin{equation}\label{GKP}
\left(-  u_{xx} +D^{-2}_{x}u_{yy}+   \tau u-f(u)\right)_{x}=0,
\,\, \mbox{in} \,\,  \mathbb{R}^N.
\end{equation}
In the pioneering work, De Bouard and Saut \cite{DS0, DS} considered equation \eqref{GKP} and treated a nonlinearity $f(u)=|u|^p u$ assuming that $p=\frac{m}{n}$, with $m$ and $n$ relatively prime, and $n$ is odd and $1 \leq p<4$, if $N=2$, or $1 \leq p<4 / 3$, if $N=3$. In \cite{DS0, DS}, De Bouard and Saut obtained existence and nonexistence of solitary waves for equation \eqref{GKP} by using constrained minimization method and concentration-compactness principle \cite{Lions}.  Willem [12] extended the results in \cite{DS0} to the case $N=2$ and $f(u)$ is a continuous function satisfying some structure conditions by using the mountain pass theorem,  then Wang and Willem \cite{WW} obtained multiple solitary waves for equation \eqref{GKP} in one-dimensional spaces via the Lyusternik-Schnirelman category theory.  In \cite{LS}, by using the variational method, Liang and Su proved the existence of nontrivial solitary waves for equation \eqref{GKP} with
$f(x, y, u)=Q(x, y)\vert u\vert^{p-2}u$ and $N\geq 2$, where $Q\in C(\mathbb{R}\times\mathbb{R}^{N-1}, \mathbb{R})$ satisfying some assumptions and $2<p<N^{*}=\frac{2(2N-1)}{2N-3}$, while Xuan \cite{Xu} dealt with the case where $f(u)$ satisfies some superlinear conditions in higher dimension. Applying the linking theorem stated in \cite{SZ}, He and Zou in \cite{HZ} studied existence of nontrivial solitary waves for the generalized Kadomtsev-Petviashvili equation in multi-dimensional spaces, similar results also see \cite{Z}. For further results about the generalized Kadomtsev-Petviashvili equation, we could see \cite{AMP, BIN, B, F, GW, HNS, IM, KP, S, SZ, TG, T, WW, XWD, ZXM} and the references therein.

Recently, Alves and Miyagaki in \cite{AM} established some results concerning the existence, regularity and concentration
phenomenon of nontrivial solitary waves for a class of generalized variable
coefficient Kadomtsev-Petviashvili equations in $\mathbb{R}^2$. After that, Figueiredo and Montenegro \cite{FM} proved the existence of multiple solitary waves for a generalized Kadomtsev-Petviashvili equation with a potential in $\mathbb{R}^2$. In \cite{FM}, the authors showed that the number of solitary waves corresponds to the number of global minimum points of the potential when a parameter is small enough.  Inspired by the papers mentioned above,  in this paper, we will study the existence and concentration of the nontrivial solitary waves for
the generalized Kadomtsev-Petviashvili equation $(P_{\varepsilon}^{*})$ in $\mathbb{R}^2$. Throughout the paper, we make the following assumptions on the potential $V$:
\begin{itemize}
	\item[\rm ($V_{1}$)] $V\in C^{2}(\mathbb{R}^2, \mathbb{R}) \cap L^{\infty}(\mathbb{R}^2, \mathbb{R})$ and $\displaystyle \inf_{(x, y)\in\mathbb{R}^2}V(x,y)=V_{0}>0$;

	\item[\rm ($V_{2}$)]There exists an open and bounded set $\Lambda\subset \mathbb{R}^{2}$ satisfying
     $$
       0<V_{1}=\min_{(x, y)\in \overline{\Lambda}}V(x,y)<\min_{(x, y)\in\partial \Lambda}V(x,y);
      $$

\item[\rm ($V_{3}$)] $V-1, V_{x}, V_{xx} \in L^{1}(\mathbb{R}^2) \cap L^{\infty}(\mathbb{R}^2)$.
\end{itemize}
Without loss of generality, we can assume that $(0,0) \in \Lambda$ and $V(0,0)=\displaystyle \inf_{(x, y)\in \overline{\Lambda}}V(x,y)$.	\\
Hereinafter, we suppose that $f:\mathbb{R} \to \mathbb{R}$ is an odd and  $C^2$ function that satisfies the following conditions:
\begin{itemize}

\item [$(f_1)$]   $f \in C^2(\R^2)$ \, with \, $f(0)=f'(0)=0$;\\
\item [$(f_2)$]   There exist $C>0$  and $p \in (1,2)$  such that
$$
|f''(t)|\leq C |t|^{p-1}, \,\,\, \mbox{for all}\,\,\, t \in \R;
$$

\item[\rm ($f_{3}$)] There exists  $\theta \in (2, 6)$ such that
$$
0<\theta F(t)=\theta \int^{t}_{0}f(r)dr \leq tf(t) \,\,\, \mbox{for all}
\,\,\, t\neq 0;
$$

\item[\rm ($f_{4}$)]
$$
t\to \frac{f(t)}{\vert t\vert} \ \  \mbox{is strictly increasing on} \ \ (-\infty, 0) \ \  \mbox{and} \ \ (0, \infty).
$$	
\end{itemize}

Now, we are ready to state our main result.
\begin{theorem}\label{gio1}
Suppose that the conditions $(V_{1})-(V_{2})$ and $(f_{1})-(f_{4})$ hold. Then, there exists $\varepsilon_{0}>0$ such that, for any $\varepsilon\in (0, \varepsilon_{0})$, problem $(P_{\varepsilon}^{*})$ has a nontrivial solution  $v_{\varepsilon}$. Moreover, if $\eta_{\varepsilon}\in \mathbb{R}^{2}$ is a global maximum point of $\vert v_{\varepsilon}\vert$  and $\varepsilon_n \to 0$ as $n \to \infty$,  then $\eta_{\varepsilon_n} \to (x_1, y_1)\in \Lambda$ with $V(x_1, y_1)=V_1$.
\end{theorem}

In order to prove Theorem \ref{gio1}, an important difficulty is lack of the compactness. This type of difficulty also occurs in the study of existence of standing wave solutions $\psi(x, t):=e^{-iEt/\hbar}u(x)$, with $E\in \mathbb{R}$, for the nonlinear Schr\"odinger  equation
 $$
 i \hbar \partial_t \Psi = - \epsilon^{2}\Delta \Psi + (W(x	)+w)\Psi- f(|\Psi|)\Psi, \ \Psi:[0, \infty) \times \mathbb{R}^{N} \rightarrow \mathbb{C}, \ N \geq 2 \eqno{(NSE)}
 $$
 that appears in some important physical applications involving quantum mechanics. In \cite{DF},  del Pino and Felmer developed a penalization method to prove existence of standing wave solutions for $(NSE)$ that permits to overcome the lack of compactness by working with a modified problem.  Currently for equation $(P_{\varepsilon})$, we will face some new difficulties, such as: \\

\noindent  1) \,\, For $u\in X_{\epsilon}$, where $X_{\epsilon}$ is the working space and its definition can be found below, it is not clear whether
$u^{+}\in X_{\epsilon}$ or  $u^{-}\in X_{\epsilon}$, where
$$
u^{+}=\max\{u(x), 0\} \quad \mbox{and} \quad u^{-}=\min\{u(x), 0\}.
$$
\noindent 2) \,\, If $\phi_{R}\in C^{\infty}(\mathbb{R}^{2})$ is a cut-off function, it is also unclear whether $\phi_{R} u \in X_{\epsilon}$.

Having the above remarks in mind, in order to apply the penalization method in \cite{DF} to study equation $(P_{\varepsilon})$, we need to adapt the penalization method in \cite{DF}, however this is not a simple thing to do, for example the reader is invited to see that the proof of Lemma \ref{compactness} is different from that explored in \cite{DF}. On the other hand, in \cite{DF}, in order to show  the solutions obtained for the modified problem are solutions of the original problem when $\epsilon>0$ is sufficiently small in the study of nonlinear Schr\"odinger  equation, the authors used the properties 1) and 2) mentioned above, then in the present paper it was necessary to develop a new approach that involves the Fourier transform of a tempered distribution to solve the problem.
In order to use the Fourier transform of a tempered distribution, it was also necessary to do a new adjustment in the deformation introduced in \cite{DF}, because in our case the deformation needs to have more regularity. Hence, the reader is invited to see the assumptions that we are considering on potential $V$ and the nonlinearity $f$ are more stronger than \cite{DF}, for example we assume that $V\in C^{2}(\mathbb{R}^2, \mathbb{R}) \cap L^{\infty}(\mathbb{R}^2, \mathbb{R})$ and $f:\mathbb{R} \to \mathbb{R}$ is a $C^2$ function with $\vert f''(t)\vert\leq C |t|^{p-1}$ for all $t \in \R$ and $p \in (1,2)$. Here we assume that $f$ is an odd function to avoid more technicality in the definition of function $g$, see \eqref{1.7} below. In fact, this condition is not necessary, for more details, we could refer \cite{AM}.  We would also like to point out that for the existence of solution for the modified problem $(G_{\varepsilon})$, see Section 2, it is enough to consider $p \in(0,4)$ in $(f_2)$, while the restriction $p \in(1, 4)$ is due to a technical difficulty to regularize these solutions, for more details see Section 3. Finally, the restriction $p \in (1,2)$ is used to prove that the solutions of the modified problem are in fact solution of the original problem when $\epsilon>0$ is small, see Section 4.

After these comments, we would like to give a natural question:
is it possible to prove  Theorem \ref{gio1} under the assumptions on $V$ and $f$ found in \cite{DF}? In our opinion, this is an interesting question.

The paper is organized as follows. In Section 2, we present our functional setting and study the modified problem, this is a key point in our approach. In Section 3, we study the regularity of solutions. Finally, in Section 4, we give the proof of Theorem \ref{gio1}.\\

\noindent \textbf{Notation:} Throughout the paper, unless explicitly stated, the symbols $C$ and $C_{i}$ will always denote a generic positive constant, which may vary from line to line. The symbols " $\rightarrow$ " and " $\rightharpoonup$ " denote, respectively, strong and weak convergence. Moreover, we denote by $|\cdot|_{q}$ the usual norm of Lebesgue space $L^{q}\left(\mathbb{R}^{2}\right)$, for $q \in[2,+\infty]$. Finally, for all $(x, y) \in \mathbb{R}^{2}$ and $r>0$, we denote by $B_{r}((x, y))$ the ball centred in $(x, y)$ with radius $r$. For the measurable set $A\subset \mathbb{R}^2$, $\vert A\vert$ denotes the Lebesgue measure of the set $A$.

\section{The variational framework} \label{var}

To begin with, we have that after a rescaling,  problem $(P_{\varepsilon}^{*})$ is equivalent to the problem below
$$
\left(-  u_{xx} +D^{-2}_{x}u_{yy}+   V(\varepsilon x,\varepsilon y)u-f(u)\right)_{x}=0,
\,\, \mbox{in} \,\,  \mathbb{R}^2. \eqno{(P_{\varepsilon})}
$$

Arguing as in Willem \cite[Chapter 7]{Willem}, in the set $Y=\{ g_x : g \in C^{\infty}_0 (\mathbb{R}^{2}) \}$ we define the inner product
\begin{eqnarray}\label{Produtointerno}
(u,v)_{\varepsilon}= \displaystyle\int_{\mathbb{R}^{2}} \left(u_x v_x + D^{-1}_{x} u_{y} D^{-1}_{x} v_{y} + V(\varepsilon x, \varepsilon y)u v\right)dx dy
\end{eqnarray}
and the corresponding norm
\begin{eqnarray}\label{Norma}
\|u\|_{\varepsilon}= \left( \displaystyle\int_{\mathbb{R}^{2}} \left(|u_{x}|^{2} + |D^{-1}_{x} u_{y}|^{2} + V(\varepsilon x, \varepsilon y)u^{2}\right)dx dy\right)^{\frac{1}{2}}.
\end{eqnarray}
We proceed to define the space $X_{\varepsilon}$ of functions that we will work with. A function $u:\mathbb{R}^{2}\to\mathbb{R}$ belongs to $X_{\varepsilon}$ if there exists $(u_n) \subset Y$ such that
$$
u_n \to u \ \  a.e. \ \ \text{in}\ \ \mathbb{R}^{2} \ \ \mbox{and} \ \ \|u_j -u_k\|_{\varepsilon}\to 0 \ \ \mbox{as} \ \ j, k \to \infty.
$$
The space $X_{\varepsilon}$ with inner product (\ref{Produtointerno}) and norm (\ref{Norma}) is a Hilbert space. By a solution of $(P_{\varepsilon})$ we mean a function $u \in X_{\varepsilon}$ such that
$$
(u,\phi)_{\varepsilon} - \displaystyle\int_{\mathbb{R}^{2}} f(u)\phi dx dy = 0  \,\,\,\, \mbox{for all} \,\,  \phi\in X_{\varepsilon}.
$$

As in \cite{DF}, to study problem $(P_{\varepsilon})$ by variational methods, we modify suitably the nonlinearity $f$ so that, for \linebreak  $\varepsilon>0$ small enough, the solutions of the modified problem are also solutions of the original problem $(P_{\varepsilon})$. However,  the modification that will introduce is different from that found in \cite{DF}, because in the present paper we need more regularity in the nonlinearity that appears in the modified problem.

First of all, without loss of generality, we assume that
$$
V(0, 0)=V_{1}=\min _{(x,y) \in \overline{\Lambda}} V(x,y).
$$
Then, choose  $k>0$ such that $k>\frac{\theta}{\theta-2}$, let $a,\alpha_0>0$ and $\delta>0$ small enough verifying $\frac{f(a)}{a}=(\frac{V_{0}}{k}-\delta)>0, \frac{f(a+2\alpha_0)}{a+2\alpha_0}=\frac{V_{0}}{k}$, where $V_{0}$ is given in $(V_{1})$. Using these numbers, we consider the function
\begin{equation*}
	\tilde{f}(t):=f(t)\phi(t)+(1-\phi(t))\frac{V_0}{k}t, \quad \forall t \in \mathbb{R},
\end{equation*}
where $\phi \in C^\infty(\mathbb{R})$ is even, nonnegative and nonincreasing  for $t \geq 0$ and  satisfying
$$
0 \leq \phi(t) \leq 1 \quad \forall t \in \mathbb{R}, \quad \phi(t)=1 \quad \mbox{for} \quad 0\leq t \leq a+\alpha_0 \quad \mbox{and} \quad \phi(t)=0 \quad \mbox{for} \quad t \geq a+2\alpha_0.
$$

From the definition of $\tilde{f}$, it follows that $\tilde{f}\in C^{2}(\mathbb{R})$,  $\tilde{f}$ is odd and $\tilde{f}(t) \leq \frac{V_0}{k}t$ for all $t \geq 0$. Hereafter,  for the above $\delta>0$ that could be decreased if necessary, we define the set
$$
\overline{\Lambda}_\delta=\{(x,y) \in \R^2\,:\, \text{dist}((x,y),\overline{\Lambda})\leq \delta\}
$$
 such that
$$
V(0,0)+\delta < \min_{\overline{\Lambda}_\delta \setminus \Lambda}V(x,y).
$$

Associated with the above set, we also consider the function $\chi \in C_0^{\infty}(\mathbb{R}^2)$ satisfying
$$
0 \leq \chi(x,y)\leq 1 \quad \forall (x,y) \in \R^2, \quad \chi(x,y)=1 \quad \forall (x,y) \in \overline{\Lambda} \quad \mbox{and} \quad \chi(x,y)=0 \quad \forall (x,y) \not\in \overline{\Lambda}_\delta.
$$
Now we introduce the penalized nonlinearity $g:\mathbb{R}^{2}\times \mathbb{R}\rightarrow \mathbb{R}$ given by
\begin{equation}
	\label{1.7}
	g(x, y, t):=\chi(x,y)f(t)+(1-\chi(x,y))\tilde{f}(t),
\end{equation}
and $\displaystyle G(x, y, t):=\int_{0}^{t}g(x, y, s)ds$. In view of  $(f_{1})$--$(f_{4})$, we deduce that $g$ is odd in $t\in \mathbb{R}$ and $C^2$ in $(t,x, y)\in \mathbb{R}^{3}$ satisfying the following properties:
\begin{itemize}
 \item[$(g_{0})$]
     $g(x,y,t)\leq \delta t+f(t)$ for any $t \geq 0$ and there exists a constant $C>0$ such that
     $$ \vert g_x(x,y,t)\vert\leq C\vert t\vert^{p+1},|g_{xx}(x,y,t)|\leq C\vert t\vert^{p+1}, |g_{tt}(x,y,t)| \leq C \vert t\vert^{p-1}, \quad \forall t \geq 0;
     $$

\item[$(g_{1})$]$\underset{t\rightarrow 0}{\lim}\frac{ g(x, y,  t)}{t}=0$ uniformly in $(x, y)\in \mathbb{R}^{2}$; \\

\item[$(g_{2})$] $0<\theta G(x, y, t)\leq  g(x, y, t)t$, for each $(x, y)\in \overline{\Lambda}$, $t\neq 0$;

\item[$(g_{3})$] $0<2G(x, y, t)\leq  g(x, y, t)t\leq V_{0}t^{2}/k$, \, for each $(x, y) \not\in \overline{\Lambda}_\delta$, $t\neq 0$;

\item[$(g_{4})$] For each $(x, y)\in \overline{\Lambda}$,  the function  $t\mapsto \frac{g(x, y, t)}{\vert t\vert}$ is strictly increasing on $(-\infty, 0)$ and $(0, \infty)$.  Moreover, for each $(x, y)\in \Big(\overline{\Lambda}_\delta\Big)^c$, the function  $t\mapsto \frac{g(x, y, t)}{\vert t\vert}$  is strictly increasing on $(-a-\alpha_0, 0)$ and $(0, a+\alpha_0)$.
 \end{itemize}
From $(g_{3})$, one has

\begin{itemize}
\item[$(g_{5})$] $L(x, y, t)=(1-\chi(x, y))\Big(\frac{1}{2}V( x,  y) t^{2}-\tilde{F}(t)\Big)\geq (1-\chi(x, y))(\frac{1}{2}-\frac{1}{2k}) V( x,  y) t^{2}\geq 0$, for each $(x, y) \in \mathbb{R}^2$, $t\in\mathbb{R}$;
\item[$(g_{6})$] $M(x, y, t)=(1-\chi(x, y))(V( x,  y) t^{2}-\tilde{f}(t)t)\geq (1-\chi(x, y))(1-\frac{1}{k}) V( x,  y)t^{2}\geq 0$, for each $(x, y) \in \mathbb{R}^2$, $t\in\mathbb{R}$.
 \end{itemize}

Now we consider the  modified problem
$$
\left(-  u_{xx} +D^{-2}_{x}u_{yy}+   V(\varepsilon x,\varepsilon y)u-g(\varepsilon x,\varepsilon y, u)\right)_{x}=0
\,\,\,\, \mbox{in} \,\,\,  \mathbb{R}^2. \eqno{(G_{\varepsilon})}
$$
Note that, if  $u_{\varepsilon}$
is a solution of problem $(G_{\varepsilon})$ with
$$
	\vert u_{\varepsilon}(x, y)\vert \leq a\quad\text{for all } (x, y)\not\in 	(\overline{\Lambda}_\delta)_{\varepsilon},
$$
where
$$
	(\overline{\Lambda}_\delta)_{\varepsilon}:=\Big\{(x, y)\in \mathbb{R}^{2}: \varepsilon (x, y)\in \overline{\Lambda}_\delta\Big\},
$$
then $u_{\varepsilon}$ is a solution of problem $(P_{\varepsilon})$ .

The energy functional associated to problem $(G_{\varepsilon})$  is
$$
J_{\varepsilon}(u)=\frac 12 \|u\|^{2}_{\varepsilon} - \displaystyle\int_{\mathbb{R}^{2}}G(\varepsilon x, \varepsilon y, u) dx dy.
$$
It is standard to prove that $J_{\varepsilon}\in C^{1}(X_{\varepsilon}, \mathbb{R})$ and its critical points are the weak solutions of the modified
problem $(G_{\varepsilon})$. We would like to point out that we can assume $p \in (0, 4)$ and $p \in (1, 4)$ in the present section and Section 3 respectively. The restriction $p \in (1,2)$ is only used in Section 4.

Now we show that the functional $J_{\varepsilon}$ satisfies the mountain pass geometry (see \cite{Willem}).

\begin{lemma}\label{geometriadamontanha}
The functional $J_{\varepsilon}$ satisfies the following properties. \\
(i) There exist $r$, $\beta> 0$ such that
$$
J_{\varepsilon}(u)\geq \beta \,\,\, \mbox{with} \,\,\,
\|u\|_\varepsilon=r.
$$
(ii) There exists $ \| e \|_\varepsilon > r$ verifying $J_{\varepsilon}(e)<0$.
\end{lemma}
\noindent\textbf{Proof.}
(i)\,\, By $(g_{1})$, $(f_{2})$ and the definition of $g(x, y, t)$, for any $\zeta>0$ small, there exists $C_{\zeta}>0$ such that
\begin{equation*}
\label{3.1}
G(\varepsilon x, \varepsilon y, u)\leq\frac{\delta +\zeta}{2} u^{2}+C_{\zeta}\vert u\vert^{p+2}\,\, \,\text{for all }\,\, (x, y)\in \mathbb{R}^{2},
\end{equation*}
where  $p \in (0, 4)$. By the  Sobolev embedding,
\begin{align*}
J_{\varepsilon}(u)
& \geq
\frac{1}{2}\Vert u_{n}\Vert^{2}_{\varepsilon}-\frac{(\delta +\zeta)}{2}\int_{\mathbb{R}^{2}}u^{2}dxdy-C_{\zeta}\int_{\mathbb{R}^{2}}\vert u\vert^{p+2}dxdy\\
&\geq
\frac{1}{4}\Vert u_{n}\Vert^{2}_{\varepsilon}-C_{1}\Vert u_{n}\Vert^{p+2}_{\varepsilon}.
\end{align*}
Hence we can choose some $\beta, r>0$ such that $J_{\varepsilon}(u)\geq \beta$ if $\Vert u\Vert_{\varepsilon}=r$ small.\\
(ii)\,\, Let $u\in X_{\varepsilon}\setminus\{0\}$ with $\text{supp}(u)\subset \Lambda_{\varepsilon}$. By the definition of $g$,
$$
J_{\varepsilon}(t u)=
\frac{t^{2}}{2} \|u\|^{2}_{\varepsilon} - \displaystyle\int_{\mathbb{R}^{2}}G(\varepsilon x, \varepsilon y, tu) dx dy = \frac{t^{2}}{2} \|u\|^{2}_{\varepsilon} - \int_{\overline{\Lambda}_{\varepsilon}} F(t u) dx dy.
$$
This together with $(f_{3})$ implies that $J_{\varepsilon}(t u)\rightarrow -\infty$ as $t\rightarrow +\infty$ and the conclusion follows. \;\;\;\; $\square$\\

The main feature of the modified functional is that it satisfies the compactness condition,
as we can see in the following results.

\begin{lemma}\label{compactness}
For any fixed $\varepsilon>0$, the functional  $J_{\varepsilon}$ verifies the $(PS)_c$ condition.
\end{lemma}

\noindent\textbf{Proof.}
Assume that $(u_{n})\subset X_{\varepsilon}$ is a $(PS)_{c}$ sequence for $J_{\varepsilon}$, that is,
$$J_{\varepsilon}(u_{n})\rightarrow c\,\,\,\text{and}\,\,\,J'_{\varepsilon}(u_{n})\rightarrow 0.$$
From this, there exists $n_{0}\in \mathbb{N}$ such that
\begin{equation} \label{INE1}
J_{\varepsilon}(u_{n})
-\frac{1}{\theta}J_{\varepsilon}'(u_{n})(u_{n})\leq c+1+o_{n}(1)\|u_n\|_\varepsilon,\,\,\, \forall n\geq n_{0}.
\end{equation}
On the other hand, by $(g_{2})$ and $(g_{3})$,
\begin{align*}
J_{\varepsilon}(u_{n})
-\frac{1}{\theta}J_{\varepsilon}'(u_{n})(u_{n})=&
\Big(\frac{1}{2}-\frac{1}{\theta}\Big)\Vert u_{n}\Vert^{2}_{\varepsilon}+ \int_{\mathbb{R}^{2}}\Big(\frac{1}{\theta}g(\varepsilon x, \varepsilon y,  u_{n})u_{n}-G(\varepsilon x, \varepsilon y,  u_{n})\Big)dxdy\\
\geq&\Big(\frac{1}{2}-\frac{1}{\theta}\Big)\Vert u_{n}\Vert^{2}_{\varepsilon}
+\int_{\mathbb{R}^2}
(1-\chi(\epsilon x,\epsilon y))\Big(\frac{1}{\theta}\tilde{f}(u_n)u_{n}-\tilde{F}(u_{n})\Big)dxdy\\
\geq&\Big(\frac{1}{2}-\frac{1}{\theta}\Big)\Vert u_{n}\Vert^{2}_{\varepsilon}-\int_{\mathbb{R}^2}(1-\chi(\epsilon x,\epsilon y))
\tilde{F}(u_{n})dxdy\\
\geq&\Big(\frac{1}{2}-\frac{1}{\theta}\Big)\Vert u_{n}\Vert^{2}_{\varepsilon}-\frac{1}{2k}\int_{\mathbb{R}^{2}}
V_{0}\vert u_{n}\vert^{2}dxdy\\
\geq&\Big(\frac{1}{2}-\frac{1}{\theta}-\frac{1}{2k}\Big)\Vert u_{n}\Vert^{2}_{\varepsilon},
\end{align*}
that is,
\begin{equation} \label{INE2}
J_{\varepsilon}(u_{n}) -\frac{1}{\theta}J_{\varepsilon}'(u_{n})(u_{n}) \geq \Big(\frac{1}{2}-\frac{1}{\theta}-\frac{1}{2k}\Big)\Vert u_{n}\Vert^{2}_{\varepsilon}.
\end{equation}	
From (\ref{INE1})-(\ref{INE2}),
$$
c+1+o_{n}(1)\|u_n\|_\varepsilon\geq\Big(\frac{1}{2}-\frac{1}{\theta}-\frac{1}{2k}\Big)\Vert u_{n}\Vert^{2}_{\varepsilon},\,\, \forall n\geq n_{0}.
$$
Since $k>\frac{\theta}{\theta-2}$,  the above inequality guarantees that $\{u_{n}\}$ is a bounded sequence in $X_{\varepsilon}$. Thus, up to a subsequence, $u_{n}\rightharpoonup u$
in $X_{\varepsilon}$, $u_{n}\rightarrow u$ in $L^{r}_{\rm loc}(\mathbb{R}^{2})$ for $1\leq r< 6$, and $u_{n}(x, y)\rightarrow u(x, y)$ a.e. in $\mathbb{R}^{2}$ as $n\rightarrow+\infty$. These combined with Lebesgue Dominated Convergence Theorem ensure that for each $\psi\in Y$, we have the limit below
$$
\displaystyle \int_{\mathbb{R}^2}g(\varepsilon x, \varepsilon y,  u_{n})\psi dx dy\rightarrow
\displaystyle \int_{\mathbb{R}^2}g(\varepsilon x, \varepsilon y,  u)\psi dx dy.
$$
The above limit together with the weak convergence yields that $J_{\varepsilon}'(u)(\psi)=0$ for all $\psi \in Y$, that is, $J_{\varepsilon}'(u)=0$.\\

Since $J_{\varepsilon}'(u)(u)=0$, one has
$$
\displaystyle\int_{\mathbb{R}^{2}} \left(u_{x}^{2} + |D^{-1}_{x} u_{y}|^{2} \right)dx dy+\displaystyle\int_{\mathbb{R}^2}\chi(\epsilon x, \epsilon y) V(\varepsilon x, \varepsilon y)|u|^{2}dx dy+\displaystyle\int_{\mathbb{R}^2}M(\varepsilon x, \varepsilon y,  u)dxdy=\int_{\mathbb{R}^2}\chi(\epsilon x, \epsilon y) f(u)udx dy,
$$
where $M$ was fixed in $(g_6)$. On the other hand, by $J_{\varepsilon}'( u_{n})(u_{n})=o_{n}(1)$,
\begin{align*}
&\int_{\mathbb{R}^{2}} \left((u_{n})_{x}^{2} + |D^{-1}_{x} (u_{n})_{y}|^{2} \right)dx dy+\displaystyle\int_{\mathbb{R}^2} \chi(\epsilon x,\epsilon y) V(\varepsilon x, \varepsilon y)|u_{n}|^{2}dx dy+\displaystyle\int_{\mathbb{R}^2}M(\varepsilon x, \varepsilon y,  u_{n})dxdy \\
=&\displaystyle \int_{\mathbb{R}^2}\chi(\epsilon x, \epsilon y)f(u_{n})u_{n} dx dy.
\end{align*}
Since the domain $(\overline{\Lambda}_\delta)_{\varepsilon}$ is bounded, the compactness of the Sobolev embedding and the Lebesgue Dominated Convergence theorem lead to
$$
\lim_{n\to +\infty}\int_{\mathbb{R}^2}\chi(\epsilon x, \epsilon y)f(u_{n})u_{n} dx dy= \int_{\mathbb{R}^2}\chi(\epsilon x, \epsilon y)f(u)udx dy
$$
and
\begin{align}\label{eq0}
\lim_{n\to +\infty}\int_{\mathbb{R}^2} \chi(\epsilon x, \epsilon y) V(\varepsilon x, \varepsilon y)\vert u_{n}\vert^{2}dx dy= \int_{\mathbb{R}^2}\chi(\epsilon x, \epsilon y)  V(\varepsilon x, \varepsilon y)u^{2}dx dy.
\end{align}
Therefore,
\begin{align*}
&\lim_{n\to +\infty}\Big(\displaystyle\int_{\mathbb{R}^{2}} \left(|(u_{n})_{x}|^{2} + |D^{-1}_{x} (u_{n})_{y}|^{2} \right)dx dy+\displaystyle\int_{\mathbb{R}^2}M(\varepsilon x, \varepsilon y,  u_{n})dxdy\Big)\\
=&\displaystyle\int_{\mathbb{R}^{2}} \left(|u_{x}|^{2} + |D^{-1}_{x} u_{y}|^{2} \right)dx dy+\displaystyle\int_{\mathbb{R}^2}M(\varepsilon x, \varepsilon y,  u)dxdy.
\end{align*}
From $(g_6)$, we know that $M(x, y, t)\geq 0$ for each $(x, y)\in \mathbb{R}^2$ and $t\in \mathbb{R}$. Then, by the Fatou's lemma,
\begin{align*}
&\liminf_{n\to +\infty}\int_{\mathbb{R}^2}M(\varepsilon x, \varepsilon y,  u_{n})dxdy\geq\displaystyle\int_{\mathbb{R}^2}M(\varepsilon x, \varepsilon y,  u)dxdy.
\end{align*}
On the other hand, since $\displaystyle \int_{\mathbb{R}^{2}} \left(|u_{x}|^{2} + |D^{-1}_{x} u_{y}|^{2} \right)dx dy$
is continuous and convex in $X_{\varepsilon}$, we must have
\begin{align*}
\liminf_{n\to +\infty}\Big(\displaystyle\int_{\mathbb{R}^{2}} \left(|(u_{n})_{x}|^{2} + |D^{-1}_{x} (u_{n})_{y}|^{2} \right)dx dy
\geq\displaystyle\int_{\mathbb{R}^{2}} \left(|u_{x}|^{2} + |D^{-1}_{x} u_{y}|^{2} \right)dx dy.
\end{align*} \\
Hence,
\begin{align}\label{eq1}
\lim_{n\to +\infty}\int_{\mathbb{R}^2}M(\varepsilon x, \varepsilon y,  u_{n})dxdy=\displaystyle\int_{\mathbb{R}^2}M(\varepsilon x, \varepsilon y,  u)dxdy
\end{align}
and
\begin{align}\label{eq2}
\lim_{n\to +\infty}\Big(\displaystyle\int_{\mathbb{R}^{2}} \left(|(u_{n})_x|^{2} + |D^{-1}_{x} (u_{n})_{y}|^{2} \right)dx dy\Big)
=\displaystyle\int_{\mathbb{R}^{2}} \left(|u_{x}|^{2} + |D^{-1}_{x} u_{y}|^{2} \right)dx dy.
\end{align}
From \eqref{eq1}, the definition of $M$ and a variant of the Lebesgue Dominated Convergence Theorem, we have
$$
\lim_{n\to +\infty}\displaystyle\int_{\mathbb{R}^2} (1-\chi(\epsilon x,\epsilon y)) V(\varepsilon x, \varepsilon y)\vert u_{n}\vert^{2}dx dy=\displaystyle\int_{\mathbb{R}^2}  (1-\chi(\epsilon x,\epsilon y)) V(\varepsilon x, \varepsilon y)u^{2}dx dy.
$$
Gathering this limit with \eqref{eq0}, we deduce that
\begin{align}\label{eq4}
\lim_{n\to +\infty}\int_{\mathbb{R}^{2}}  V(\varepsilon x, \varepsilon y)|u_{n}|^{2}dx dy= \int_{\mathbb{R}^{2}}  V(\varepsilon x, \varepsilon y)u^{2}dx dy.
\end{align}
From \eqref{eq2} and \eqref{eq4}, it yields
$$
\lim_{n\to +\infty}\Vert u_{n}\Vert^{2}_{\varepsilon}=\Vert u\Vert^{2}_{\varepsilon}.
$$
Since $X_{\varepsilon}$ is a Hilbert space and $u_{n}\rightharpoonup u$
in $X_{\varepsilon}$, the above limit yields
$$
u_{n}\rightarrow u\,\,\, \text{in}\,\,\,X_{\varepsilon},
$$
proving the lemma. \;\;\; $\Box$

\medskip

\begin{theorem}
 The functional $J_{\varepsilon}$ has a nontrivial critical point $u_{\varepsilon} \in X_{\varepsilon}$ such that
\begin{align}\label{equ5}
J_{\varepsilon}\left(u_{\varepsilon}\right)=c_{\varepsilon} \quad \text { and } \quad J_{\varepsilon}^{\prime}\left(u_{\varepsilon}\right)=0,
\end{align}
where $c_{\varepsilon}$ denotes the mountain pass level associated with $J_{\varepsilon},$ that is,
\begin{equation*}\label{d0}
	c_\varepsilon=\displaystyle \inf_{\eta \in \Gamma}\displaystyle \max_{t\in [0,1]} \displaystyle  J_{\varepsilon}(\eta(t)),
\end{equation*}
where $\Gamma =\{\eta \in C([0,1], X_{\varepsilon}): \eta(0)=0 \ \mbox{and} \ J_{\varepsilon}(\eta(1))<0\}$.
\end{theorem}
\noindent\textbf{Proof.} Lemmas \ref{geometriadamontanha} and  \ref{compactness} permit to apply the Mountain Pass Theorem due to Ambrosetti and Rabinowitz (\cite[ Theorem 1.15]{Willem}) to conclude that $c_\varepsilon$ is a critical level for $J_{\varepsilon}$. \;\;\; $\Box$

Related to the case $\varepsilon=0$, see \cite[Proposition 2.5]{FM}, we know that there exists a nontrivial function $w \in X$ such that
$$
J_{0}(w)=c_{0} \quad \text { and } \quad J_{0}^{\prime}(w)=0
$$
where
$$
J_{0}(u) =\frac{1}{2}\int_{\mathbb{R}^{2}} \left(|u_{x}|^{2} + |D^{-1}_{x} u_{y}|^{2} + V(0, 0)u^{2}\right)dx dy-\int_{\mathbb{R}^2}F(u)dx dy
$$
and
$$
c_{0}=\inf_{u\in \mathcal{M}_{0}}J_{0}(u)=\inf _{\eta \in \Gamma} \max _{t \in[0,1]} J_{0}(\eta(t)),
$$
where $\mathcal{M}_{0}=\left\{u \in X \backslash\{0\}: J_{0}'(u)(u)=0\right\}$ and $\Gamma=\left\{\eta \in C([0,1], X): \eta(0)=0\right.$ and $\left.J_{0}(\eta(1))<0\right\}$. Hereafter $X$ is defined as the space $X_\varepsilon$ but endowed with the following norm
$$
\|u\|=  \left(\displaystyle\int_{\mathbb{R}^{2}} \left(|u_{x}|^{2} + |D^{-1}_{x} u_{y}|^{2} + V(0, 0)|u|^{2}\right)dx dy\right)^{\frac{1}{2}}.
$$

The next lemma establishes an important relation between $c_{\varepsilon}$ and $c_{0}$.

\begin{lemma} \label{lemma1}
The numbers $c_{\varepsilon}$ and $c_{0}$ verify the following inequality
$$\displaystyle  \limsup_{\varepsilon\rightarrow 0}c_{\varepsilon}\leq c_{0}.$$
\end{lemma}
\noindent\textbf{Proof.}
\noindent By \cite{FM}, we know that there exists a nontrivial function $w\in \mathcal{M}_{0}$ such that  $c_{0}=J_{0}(w)$ and $J'_{0}(w)=0$.  As $w  \not \equiv 0$, there exists $R>0$ such that
$ |\left\{ (x, y) \in B_R((0, 0))\,:\,|w(x, y)|\not= 0\right\}|>0$. Hence, for $\varepsilon>0$ small enough we have that $B_{\varepsilon R}((0, 0)) \subset \overline{\Lambda}$, and so,
$$
\displaystyle\int_{\mathbb{R}^{2}}G(\varepsilon x, \varepsilon y, t_{\varepsilon} w) dx dy \geq \displaystyle\int_{B_R((0, 0))}G(\varepsilon x, \varepsilon y, t_{\varepsilon} w) dx dy=\displaystyle\int_{B_R((0, 0))}F(t_{\varepsilon} w) dx dy,\, \,\text{for all }\,\,  t \geq 0.
$$
The above inequality combines with $(f_3)$ to give
$$
\lim_{t \to +\infty}J_\varepsilon(tw)=-\infty.
$$
From this,  there exists $t_\varepsilon>0$ such that
$$
c_\varepsilon \leq \max_{t \geq 0}J_\varepsilon(tw)=J_\varepsilon(t_\varepsilon w),
$$	
that is
$$
c_\varepsilon \leq \frac{t_{\varepsilon}^2}{2}\displaystyle\int_{\mathbb{R}^{2}} \left( |w_{x}|^{2} + |D^{-1}_{x} w_{y}|^{2} + V(\varepsilon x,\varepsilon y)|w|^{2}\right)dx dy-\displaystyle\int_{\mathbb{R}^{2}}G(\varepsilon x, \varepsilon y, t_{\varepsilon} w) dx dy.
$$	
It's easy to show that $(t_\varepsilon)$ is bounded for $\varepsilon>0$ small enough,  going if necessary to a subsequence $t_\varepsilon \to t_0>0$ as $\varepsilon \to 0$. Now, taking the limit of $\varepsilon \to 0$ and using the fact that $V\in L^{\infty}(\mathbb{R}^N, \mathbb{R}) \cap C(\mathbb{R}^N, \mathbb{R})$, the Lebesgue Dominated Convergence Theorem ensures that
$$
\limsup_{\varepsilon \to 0}c_\varepsilon \leq \frac{t_{0}^2}{2}\displaystyle\int_{\mathbb{R}^{2}} \left( |w_{x}|^{2} + |D^{-1}_{x} w_{y}|^{2} + V(0,0)|w|^{2}\right)dx dy-\displaystyle\int_{\mathbb{R}^{2}}G(0,0,t_{0} w) dx dy.
$$	
Since $G(0,0, t_0w)=F(t_0w)$, it follows that
$$
\limsup_{\varepsilon \to 0}c_\varepsilon \leq J_0(t_0w) \leq J_0(w)=c_0.
$$
\hspace{19 cm}$\Box$
\begin{lemma} \label{lemma2}
 There are $r, \gamma, \epsilon^{*}>0$ and $\left(x_{\epsilon}, y_{\epsilon}\right) \subset \mathbb{R}^{2}$ such that
$$
\int_{B_{r}\left((x_{\epsilon}, y_{\epsilon})\right)}\left|u_{\varepsilon}(x, y)\right|^{2} d x dy\geq \gamma, \,\, \text{for any}\,\, \epsilon \in\left(0, \epsilon^{*}\right) .
$$
\end{lemma}
\noindent\textbf{Proof.}
 First of all, since $\left(u_{\epsilon}\right)$ satisfies \eqref{equ5}, then from Lemma \ref{geometriadamontanha}(i), there exists $\beta>0$, which is independent of $\epsilon$, such that
\begin{align}\label{equ6}
\left\|u_{\epsilon}\right\|_{\epsilon}^{2} \geq \beta>0, \, \,\text{for all }\,\, \epsilon>0.
\end{align}
To prove the lemma, it is enough to see that for any sequence $\left(\epsilon_{n}\right) \subset(0,+\infty)$ with $\epsilon_{n} \rightarrow 0$ as $n\to \infty$, the limit below
$$
\lim _{n \rightarrow+\infty} \sup _{(x, y) \in \mathbb{R}^{2}} \int_{B_{r}((x, y))}\left|u_{\epsilon_{n}}(x, y)\right|^{2} d xdy=0
$$
does not hold for any $r>0$. Otherwise, if it holds for some $r>0$, from a Lions-type result for $X_{\epsilon}$ found in \cite[Lemma 7.4]{Willem}, one gets
$$
u_{\epsilon_{n}}(x, y) \rightarrow 0 \,\, \text {in } L^{q}\left(\mathbb{R}^{2}\right), \,\,\text {as}\,\, n\rightarrow\infty\,\,\text{for any} \,\,q \in\left(2, 6\right) .
$$
Using the above limits and \eqref{equ5}, it is easy to show that
$$
\left\|u_{\epsilon_{n}}\right\|_{\epsilon_{n}}^{2} \rightarrow 0 \quad \text { as } n \rightarrow+\infty
$$
which contradicts \eqref{equ6}. \;\;\; $\Box$\\

The next compactness result is fundamental for showing that the solutions of the modified
problem are solutions of the original problem.

\begin{lemma} \label{lea} Let $\epsilon_{n}\rightarrow 0$ and  $u_{n}\in X_{\epsilon_{n}}$ such that $J_{\epsilon_{n}}(u_{n})=c_{\epsilon_{n}}$
and $J'_{\epsilon_{n}}(u_{n})=0$. Then there exists a sequence $(x_{n}, y_{n})\subset \mathbb{R}^2$
such that $\psi_{n}(x, y)=u_{n}(x_{n}+x, y_{n}+y)$ has a convergent subsequence in $X$.
Moreover, there exists $(x_{0}, y_{0})\in\Lambda$ such that
$$
\lim_{n\rightarrow\infty}\epsilon_{n}(x_{n}, y_{n})=(x_{0}, y_{0})\quad \text{and}\quad V(x_{0}, y_{0})=V_{1}.
$$
\end{lemma}

\noindent\textbf{Proof.} From the proof of Lemma \ref{compactness}, it is
easy to see that $(\Vert u_{n}\Vert^{2}_{\varepsilon_n})$ is bounded in $\R$, from where it follows that $(u_{n})$ is also bounded in $X$. Using \cite[Lemma 7.4]{Willem},
there exist $r, \gamma>0$ and a sequence $((x_{n}, y_{n}))\subset \mathbb{R}^2$ such that
\begin{equation}\label{lion1}
\underset{n\rightarrow\infty}{\lim\sup}\int_{B_{r}((x_{n}, y_{n}))}\vert u_{n}(x, y)\vert^{2}dxdy\geq\gamma.
\end{equation}
Otherwise, we can conclude that
$$
u_{n}\rightarrow 0 \quad\text{in}\,\, L^{q}(\mathbb{R}^2),\,\, \text{for all}\,\, 2< q<6
$$
and so,
$$
\lim_{n \to +\infty}\int_{\mathbb{R}^{2}}G(\varepsilon x, \varepsilon y,  u_{n})dxdy= \lim_{n \to +\infty}\int_{\mathbb{R}^{2}}g(\varepsilon x, \varepsilon y,  u_{n})u_ndxdy=0.
$$
The above limit implies that $\|u_n\|_{\varepsilon_n}\to 0$, hence  $J_{\epsilon_n}(u_{n})=c_{\epsilon_n}\rightarrow 0$ as $\rightarrow\infty$, which is a contradiction, since $c_{\epsilon_n} \geq \beta>0$. Setting $\psi_{n}(x, y)=u_{n}(x+x_{n}, y+y_{n})$ , then
 there exists  $\psi\in X\backslash \{0\}$
such that
\begin{equation}\label{weak}
\psi_{n}\rightharpoonup \psi\quad \text{in}\,\, X
\end{equation}
and
\begin{equation}\label{positive}
\int_{B_{r}((0, 0))}\vert \psi(x, y)\vert^{2}dxdy\geq\gamma.
\end{equation}
In the sequel we will prove that the sequence $(\vert \epsilon_{n}(x_{n}, y_{n})\vert)$ is bounded in $\mathbb{R}$. To this end, it is enough to show the following claim.

\begin{claim} \label{A1} $\underset{n\rightarrow\infty}{\lim}\text{dist}\big(\epsilon_{n}(x_{n}, y_{n}), \overline{\Lambda}_\delta \big)=0$.\\

\end{claim}
Indeed, if the claim does not hold, there exist $\tau>0$ and a subsequence of $((\epsilon_{n}x_{n},\epsilon_{n}y_{n}))$, still denoted by itself, such that
$$
\text{dist}\big((\epsilon_{n}x_{n},\epsilon_{n}y_{n}),  \overline{\Lambda}_\delta \big)\geq \tau, \,\, \text{for all}\,\, n\in \mathbb{N}.
$$		
Consequently, there exists $r>0$ such that
$$
B_{r}((\epsilon_{n}x_{n}, \epsilon_{n}y_{n}))\subset (\overline{\Lambda}_\delta)^{c}, \,\, \text{for all}\,\, n\in \mathbb{N}.
$$
According to the above definition, we have
\begin{align*}
\int_{B_{\frac{r}{\epsilon_{n}}}((0, 0))} g\big(\epsilon_{n} (x_{n}, y_{n}) +\epsilon_{n}(x, y), \psi_{n}\big)\omega_{j}dxdy=\int_{B_{\frac{r}{\epsilon_{n}}}((0, 0))} \tilde{f}(\psi_{n})\omega_{j}dxdy.
\end{align*}
Using the fact that $\psi$ is a  function in $X$, there exists a sequence of  functions $(\omega_{j})\subset Y$
such that $\omega_{j}$ has a compact support in $\mathbb{R}^2$ and $\omega_{j}\rightarrow \psi$ in $X$ as $j\rightarrow \infty$. Now,
fixing $j>0$ and using $w_{j}$ as a test function, we have
\begin{align}\label{test}
&\displaystyle\int_{\mathbb{R}^{2}} \Big((\psi_{n})_x (\omega_{j})_x + D^{-1}_{x} (\psi_{n})_{y} D^{-1}_{x} (\omega_{j})_{y} + V(\epsilon_{n} (x_{n}, y_{n}) +\epsilon_{n}(x, y))\psi_{n} \omega_{j}\Big)dx dy\nonumber\\
=&\displaystyle\int_{\mathbb{R}^{2}}  g\big(\epsilon_{n} (x_{n}, y_{n}) +\epsilon_{n}(x, y), \psi_{n}\big)\omega_{j}dxdy.
\end{align}
Note that
\begin{align*}
&\displaystyle\int_{\mathbb{R}^{2}}  g\big(\epsilon_{n} (x_{n}, y_{n}) +\epsilon_{n}(x, y), \psi_{n}\big)\omega_{j}dxdy\\
=&\int_{B_{\frac{r}{\epsilon_{n}}}((0, 0))} g\big(\epsilon_{n} (x_{n}, y_{n}) +\epsilon_{n}(x, y), \psi_{n}\big)\omega_{j}dxdy+ \int_{\mathbb{R}^2\backslash B_{\frac{r}{\epsilon_{n}}}((0, 0))} g\big(\epsilon_{n} (x_{n}, y_{n}) +\epsilon_{n}(x, y), \psi_{n}\big)\omega_{j}dxdy
\end{align*}
and so,
$$
\begin{array}{ll}
&\displaystyle\int_{\mathbb{R}^{2}}  g\big(\epsilon_{n} (x_{n}+y_{n}) +\epsilon_{n}(x+y), \psi_{n}\big)\omega_{j}dxdy\\
\leq &\frac{V_{0}}{k}\displaystyle\int_{B_{\frac{r}{\epsilon_{n}}}((0, 0))}\vert \psi_{n}\vert \vert \omega_{j}\vert dxdy+\delta \displaystyle \int_{\mathbb{R}^2\backslash B_{\frac{r}{\epsilon_{n}}}((0, 0)) }|\psi_{n}|\omega_{j}|dxdy+ \int_{\mathbb{R}^2\backslash B_{\frac{r}{\epsilon_{n}}}((0, 0))} |f(\psi_{n})||\omega_{j}|dxdy.
\end{array}
$$
Therefore,
\begin{eqnarray*}
&&\displaystyle\int_{\mathbb{R}^{2}} \Big((\psi_{n})_x (\omega_{j})_x + D^{-1}_{x} (\psi_{n})_{y} D^{-1}_{x} (\omega_{j})_{y} + V(\epsilon_{n} (x_{n}+y_{n}) +\epsilon_{n}(x+y))(\psi_{n} \omega_{j})^{+}\Big)dx dy\\
&\leq &\|V\|_\infty\int_{\mathbb{R}^2}(\psi_{n} \omega_{j})^{-} dxdy+\frac{V_{0}}{k}\int_{B_{\frac{r}{\epsilon_{n}}}((0, 0))}\vert \psi_{n}\vert \vert \omega_{j}\vert dxdy \\
&&+ \delta \displaystyle \int_{\mathbb{R}^2\backslash B_{\frac{r}{\epsilon_{n}}}((0, 0))} |\psi_{n}|\omega_{j}|dxdy+ \int_{\mathbb{R}^2\backslash B_{\frac{r}{\epsilon_{n}}}((0, 0))} |f(\psi_{n})||\omega_{j}|dxdy,
\end{eqnarray*}
where $t^{+}=\max\{t,0\}$ and $t^{-}=\max\{-t,0\}$.
As $\omega_{j}$ has a compact support in $\mathbb{R}^2$ and $\epsilon_{n}\rightarrow 0$, the boundedness of $(\psi_{n})$ in $X$ ensures that
$$
\int_{\mathbb{R}^2}(\psi_{n} \omega_{j})^{-} dxdy \to \int_{\mathbb{R}^2}(\psi\omega_{j})^{-} dxdy, \;\;\; \int_{B_{\frac{r}{\epsilon_{n}}}((0, 0))}\vert \psi_{n}\vert \vert \omega_{j}\vert dxdy \to \int_{\mathbb{R}^2}\vert \psi \vert \vert \omega_{j}\vert dxdy,\,\,\text{as}\,\, n\rightarrow\infty,
$$
$$
\int_{\mathbb{R}^2\backslash B_{\frac{r}{\epsilon_{n}}}((0, 0))} |f(\psi_{n})||\omega_{j}|dxdy\rightarrow 0,\,\,\text{as}\,\, n\rightarrow\infty,
$$
$$
\int_{\mathbb{R}^2\backslash B_{\frac{r}{\epsilon_{n}}}((0, 0))} |\psi_{n}||\omega_{j}|dxdy\rightarrow 0,\,\,\text{as}\,\, n\rightarrow\infty,
$$
and
$$
 \liminf_{n \to +\infty }\displaystyle\int_{\mathbb{R}^{2}}V(\epsilon_{n} (x_{n}, y_{n}) +\epsilon_{n}(x, y))(\psi_{n} \omega_{j})^{+}dx dy \geq A\displaystyle\int_{\mathbb{R}^{2}}(\psi \omega_{j})^{+}dx dy
$$
where $A=\displaystyle \liminf_{|(x,y)|\to +\infty}V(x,y)$ if $(\epsilon_n(x_n,y_n))$ is unbounded or $A=V(x_0,y_0)$ if $\epsilon_n(x_n,y_n) \to (x_0, y_0)$ for some subsequence, and so, $A \geq V_0>0$. From this,
\begin{align*}
&\liminf_{n \to +\infty }\displaystyle\int_{\mathbb{R}^{2}} \Big((\psi_{n})_x (\omega_{j})_x + D^{-1}_{x} (\psi_{n})_{y} D^{-1}_{x} (\omega_{j})_{y} + V(\epsilon_{n} (x_{n}, y_{n}) +\epsilon_{n}(x, y))(\psi_{n} \omega_{j})^{+}\Big)dx dy\\
\geq & \int_{\mathbb{R}^{2}} \Big(\psi_x (\omega_{j})_x + D^{-1}_{x} \psi_{y} D^{-1}_{x} (\omega_{j})_{y} + A(\psi \omega_{j})^{+}\Big)dx dy, \quad\text{as}\,\, n\rightarrow\infty.
\end{align*}
Hence,
 \begin{equation*}
\displaystyle\int_{\mathbb{R}^{2}} \left(\psi_x (\omega_{j})_x + D^{-1}_{x} \psi_{y} D^{-1}_{x} (\omega_{j})_{y} + A(\psi \omega_{j})^{+}\right)dx dy\leq  \|V\|_\infty\int_{\mathbb{R}^2}(\psi\omega_{j})^{-} dxdy + \frac{V_{0}}{k}\int_{\mathbb{R}^2}\vert \psi \vert \vert \omega_{j}\vert dxdy.
\end{equation*}
Since $j$ is arbitrary, taking the limit of $j\rightarrow +\infty$, the Lebesgue Dominated Convergence theorem leads to
$$
\int_{\mathbb{R}^2}(\psi\omega_{j})^{-} dxdy \to \int_{\mathbb{R}^2}(\psi^2)^{-} dxdy=0
$$
and
$$
\int_{\mathbb{R}^2}(\psi\omega_{j})^{+} dxdy \to \int_{\mathbb{R}^2}(\psi^2)^{+} dxdy=\int_{\mathbb{R}^2}|\psi|^2 dxdy,
$$
and so,
 \begin{eqnarray*}
 \displaystyle\int_{\mathbb{R}^{2}} \left(|\psi_{x}|^{2} + |D^{-1}_{x} \psi_{y}|^{2} + A|\psi|^{2}\right)dx dy \leq \frac{V_{0}}{k}\int_{\mathbb{R}^2}|\psi|^2  dxdy,
\end{eqnarray*}
that is,
 \begin{eqnarray*}
	\displaystyle\int_{\mathbb{R}^{2}} \left(|\psi_{x}|^{2} + |D^{-1}_{x} \psi_{y}|^{2} + |\psi|^{2}\right)dx dy =0,
\end{eqnarray*}
which contradicts (\ref{positive}). This proves Claim \ref{A1}.\\

From Claim \ref{A1}, there exists a subsequence of $((\epsilon_{n}x_{n}, \epsilon_{n}y_{n}))$ and $(x_{0}, y_{0})\in  \overline{\Lambda}_\delta$  such that
 $$
 \lim_{n\rightarrow\infty} (\epsilon_{n}x_{n}, \epsilon_{n}y_{n})=(x_{0}, y_{0}).
 $$
\begin{claim} \label{A2} $(x_0,y_{0})\in \Lambda$.
\end{claim}
Indeed, by (\ref{test}),
$$
\begin{aligned}
&\int_{\mathbb{R}^{2}}\Big ((\psi_{n})_x (\omega_{j})_x + D^{-1}_{x} (\psi_{n})_{y} D^{-1}_{x} (\omega_{j})_{y} + V(\epsilon_{n} (x_{n}+y_{n}) +\epsilon_{n}(x+y))\psi_{n} \omega_{j}\Big)dx dy\\
=&\int_{\mathbb{R}^{2}}  g\big(\epsilon_{n} (x_{n}, y_{n}) +\epsilon_{n}(x, y), \psi_{n}\big)\omega_{j}dxdy.\\
\end{aligned}
$$
From the Lebesgue Dominated Convergence theorem, one has
$$
\int_{\mathbb{R}^{2}}  g(\epsilon_{n} (x_{n}, y_{n}) +\epsilon_{n}(x, y), \psi_{n})\omega_{j}dxdy \to \int_{\mathbb{R}^{2}}  {g}\big(x_0,y_0, \psi \big)\omega_{j}dxdy
$$
and $g(x_0,y_0,\psi)=\chi(x_0,y_0)f(\psi)+(1-\chi(x_0,y_0))\tilde{f}(\psi)$. Arguing as in the previous case, we take the limit of $n \to +\infty$ to get
$$
\displaystyle\int_{\mathbb{R}^{2}}\Big (\psi_x (\omega_{j})_x + D^{-1}_{x} \psi_{y} D^{-1}_{x} (\omega_{j})_{y} + V(x_0, y_0)\psi \omega_{j}\Big)dx dy= \int_{\mathbb{R}^{2}} g(x_0,y_0,\psi) \omega_{j}dxdy.
$$
Then, letting $j \rightarrow+\infty$, we get
$$
\displaystyle\int_{\mathbb{R}^{2}} \left(|\psi_{x}|^{2} + |D^{-1}_{x} \psi_{y}|^{2} + V(x_0, y_0)|\psi|^{2}\right)dx dy =\int_{\mathbb{R}^{2}} g(x_0,y_0,\psi) \psi dxdy.
$$
Hence,
$$
\psi \in \mathcal{N}_{V\left(x_{0}, y_{0}\right)}:=\{u\in X\backslash\{0\}: \tilde{J}^{\prime}_{V\left(x_{0}, y_{0}\right)}(u)u=0\},
$$
where $\tilde{J}_{V\left(x_{0},y_0\right)}: X \rightarrow \mathbb{R}$ is given by
$$
\tilde{J}_{V\left(x_{0}, y_{0}\right)}(u)=\frac{1}{2} \displaystyle\int_{\mathbb{R}^{2}} \left(|u_{x}|^{2} + |D^{-1}_{x} u_{y}|^{2} + V(x_0, y_0)|u|^{2}\right)dx dy-\int_{\mathbb{R}^{2}} G(x_0,y_0,u) dxdy.
$$
We claim that  $\frac{g(x,y, t)}{|t|}$ is nondecreasing for $t\not=0$ for each $(x,y) \in \mathbb{R}^2$. Since $g(x_0,y_0,\cdot)$ and $f$ are odd in $t\in\mathbb{R}$ and $\frac{f(t)}{|t|}$ is increasing for $t\not=0$, we only need to show that
$$
\left(	\frac{\tilde{f}(t)}{t} \right)'=\left(	\frac{{f}(t)}{t} \right)'\phi(t)+\left( 	\frac{{f}(t)}{t} - \frac{V_0}{k}\right)\phi'(t)\geq 0,\,\,\text{for all}\,\, t\geq 0.
$$	
Now, recall that $f(t)/t$ is increasing and $\phi(t) \geq 0$, we must have
$$
\left(	\frac{\tilde{f}(t)}{t} \right)' \geq \left( \frac{{f}(t)}{t} - \frac{V_0}{k}\right)\phi'(t).
$$	
Moreover, recall that $\phi'(t) < 0$ for $t \in (a+\alpha_0,a+2\alpha_0)$ and
$$
\frac{{f}(t)}{t} \leq \frac{{f}(a+2\alpha_0)}{a+2\alpha_0}=\frac{V_0}{k}, \,\,\text{for all}\,\, t \in (a+\alpha_0,a+2\alpha_0).
$$
From this
$$
\left( \frac{{f}(t)}{t} - \frac{V_0}{k}\right)\phi'(t) \geq 0, \,\,\text{for all}\,\, t \in (a+\alpha_0,a+2\alpha_0).
$$	
Hence,
$$
\left(	\frac{\tilde{f}(t)}{t} \right)' \geq 0, \,\,\text{for all}\,\, t \geq 0,
$$	
which showing that $\frac{\tilde{f}(t)}{t}$ is nondecreasing for $t \geq 0$. Moreover, $\frac{g(x,y, t)}{|t|}$ is increasing for $t\not= 0$ for each $(x,y) \in \mathbb{R}^2$ with $\chi(x,y)>0$. Have this in mind, it is easy to check that
$$
\tilde{c}_{V\left(x_{0}, y_{0}\right)}=\inf_{u \in  \mathcal{N}_{V\left(x_{0}, y_{0}\right)}}\tilde{J}_{V\left(x_{0}, y_{0}\right)}\left(u\right) \leq \tilde{J}_{V\left(x_{0}, y_{0}\right)}\left(\psi\right),
$$
where $\tilde{c}_{V\left(x_{0}, y_{0}\right)}$ denotes the mountain pass level associated with $\tilde{J}_{V\left(x_{0}, y_{0}\right)}$. By the definition of $g$, $(g_{2})$, $(g_{3})$, $(f_{3})$ and $(f_{4})$, we must have
$$
\begin{aligned}
&\liminf _{n \rightarrow+\infty}\Big(J_{\epsilon_{n}}(u_{n})-\frac{1}{2}J'_{\epsilon_{n}}(u_{n})(u_{n})\Big)\\
  = &\liminf _{n \rightarrow+\infty}\Big\{\int_{\mathbb{R}^{2}} \Big(\frac{1}{2} g\big(\epsilon_{n} (x_{n}, y_{n}) +\epsilon_{n}(x, y), \psi_{n}\big)\psi_{n}- G\big(\epsilon_{n} (x_{n}, y_{n}) +\epsilon_{n}(x, y), \psi_{n}\big)\Big)dxdy\Big\}\\
\geq &\int_{\mathbb{R}^{2}}\Big(\frac{1}{2} g(x_0,y_0,\psi)\psi- G(x_0,y_0,\psi)\Big)dxdy\\
=&\tilde{J}_{V\left(x_{0}, y_{0}\right)}\left(\psi\right).\\
\end{aligned}
$$
Hence,
$$
\tilde{c}_{V\left(x_{0}, y_{0}\right)} \leq \tilde{J}_{V\left(x_{0}, y_{0}\right)}\left(\psi\right) \leq \liminf _{n \rightarrow+\infty} J_{\epsilon_{n}}(u_{n})=\liminf _{n \rightarrow+\infty} c_{\epsilon_{n}}\leq c_{0}=c_{V(0, 0)} .
$$
From the definition of $g$ again, we know that $G(x_0,y_0,t) \leq \frac{\delta}{2}t^2+F(t)$, from where it follows that \linebreak $\tilde{c}_{V\left(x_{0}, y_{0}\right)} \geq {c}_{(V\left(x_{0}, y_{0}\right)-\delta)}$, where ${c}_{(V\left(x_{0}, y_{0}\right)-\delta)}$ is the mountain pass level of the functional
${J}_{(V\left(x_{0}, y_{0}\right)-\delta)}: X \rightarrow \mathbb{R}$ is given by
$$
{J}_{V\left(x_{0}, y_{0}\right)}(u)=\frac{1}{2} \displaystyle\int_{\mathbb{R}^{2}}\Big(|u_{x}|^{2} + |D^{-1}_{x} u_{y}|^{2} + {(V(x_0, y_0)-\delta)}|u|^{2}\Big)dx dy-\int_{\mathbb{R}^{2}} F(u) dxdy.
$$
This together with the inequality above yields that
$$
{c}_{(V\left(x_{0}, y_{0}\right)-\delta)} \leq c_{V(0, 0)}.
$$
Since $V(x_0, y_0)-\delta>V(0,0)$ if $(x_0,y_0) \in \overline{\Lambda}_\delta \setminus \Lambda$, we must have
$$
c_{V\left(x_{0}, y_{0}\right)-\delta} > c_{V(0,  0)},
$$
which is absurd. Thereby, $(x_0,y_0) \in \overline{\Lambda}$ and
$$
c_{V\left(x_{0}, y_{0}\right)}  \leq c_{V(0,  0)},
$$
from where it follows that
$$
V\left(x_{0}, y_{0}\right) \leq V(0, 0) \equiv V_{0}.
$$
As $V_{0} = \displaystyle \min_{\overline{\Lambda}_\delta} V(x, y)$, one gets
$$
V\left(x_{0}, y_{0}\right)=V_{0} .
$$
Moreover, recalling that $V_{0} < \min_{\overline{\Lambda}_\delta \setminus \Lambda}V(x,y),$ we conclude that $(x_{0}, y_{0}) \in \Lambda$, finishing the proof of Claim \ref{A2}.\\
	
	Now, we will prove that
	\begin{equation}\label{strong}
		\psi_{n}\rightarrow \psi\quad \text{in}\,\, X.
	\end{equation}
	First of all, we set the functions
	$$
	\tilde{\chi}_{n}^{1}(x, y)= \chi(\epsilon_n(x,y)+\epsilon_n(x_n,y_n)),
	$$
	$$
	\tilde{\chi}_{n}^{2}(x, y)= 1-\tilde{\chi}_{n}^{1}(x, y),
	$$
	$$
	h_{n}^{1}(x, y)=\left(\frac{1}{2}-\frac{1}{\theta}\right) V(\epsilon_{n} (x_{n}, y_{n}) +\epsilon_{n}(x, y)) |\psi_{n}(x, y)|^2 \tilde{\chi}_{n}^{1}(x, y) \geq 0,\,\,\text{for all}\,\, (x, y) \in \mathbb{R}^{2},
	$$
	$$
	h^{1}(x, y)=\left(\frac{1}{2}-\frac{1}{\theta}\right) V\left(x_{0}, y_{0}\right) |\psi(x, y)|^2,\,\,\text{for all}\,\, (x, y)\in \mathbb{R}^{2},
	$$
	$$
	\begin{aligned}
		h_{n}^{2}(x,y)=&\left(\left(\frac{1}{2}-\frac{1}{\theta}\right) V(\epsilon_{n} (x_{n}, y_{n}) +\epsilon_{n}(x, y)) |\psi_{n}|^{2}+\frac{1}{\theta} \tilde{f}(\psi_n)\psi_n- \tilde{F}(\psi_{n}\big)\right) \tilde{\chi}_{n}^{2}(x,y)\\
		\geq&\left(\left(\frac{1}{2}-\frac{1}{\theta}\right)-\frac{1}{2k}\right)  V(\epsilon_{n} (x_{n}, y_{n}) +\epsilon_{n}(x, y)) |\psi_{n}(x,y)|^2 \tilde{\chi}_{n}^{2}(x,y) \geq 0,
	\end{aligned}
	$$
	$$
	h_{n}^{3}(x, y)=\left(\frac{1}{\theta} f\left(\psi_{n}\right) \psi_{n}-F\left(\psi_{n}\right)\right) \tilde{\chi}_{n}^{1}(x, y) \geq 0,
	$$
	and
	$$
	h^{3}(x, y)=\frac{1}{\theta} f(\psi) \psi-F(\psi), \,\,\text{for all}\,\, (x, y) \in \mathbb{R}^{2}.
	$$
	Using the fact that
	$$
	\psi_{n}(x, y) \rightarrow \psi(x, y) \,\,\text { a.e. in } \,\,\mathbb{R}^{2} \,\,\text { and }\,\,  \epsilon_{n} (x_{n}, y_{n}) \rightarrow (x_{0}, y_{0})  \in \Lambda,\,\, \text{as}\,\, n\rightarrow\infty,
	$$
	it is easy to see that the limits below hold
	$$
	\tilde{\chi}_{n}^{1}(x, y) \rightarrow 1, \quad h_{n}^{1}(x, y) \rightarrow h^{1}(x, y), \quad h_{n}^{2}(x,y) \rightarrow 0\,\, \text { and } \,\,h_{n}^{3}(x, y) \rightarrow h^{3}(x, y) \,\,\text { a.e. in } \mathbb{R}^{2}.
	$$
	Now, by a direct computation,
	$$
	\begin{aligned}
		c_{0}  \geq \limsup _{n \rightarrow+\infty} c_{\epsilon_{n}}&=\limsup _{n \rightarrow+\infty}\left(J_{\epsilon_{n}}(u_{n})-\frac{1}{\theta} J'_{\epsilon_{n}}(u_{n})(u_{n}) \right) \\
		& = \limsup _{n \rightarrow+\infty}\Big( \int_{\mathbb{R}^{2}} \left(\frac{1}{2}-\frac{1}{\theta}\right)\left(|(\psi_{n})_{x}|^{2} + |D^{-1}_{x} (\psi_{n})_{y}|^{2} \right) dxdy+\int_{\mathbb{R}^{2}}\left(h_{n}^{1}+h_{n}^{2}+h_{n}^{3}\right) d xdy \Big)\\
		& \geq \liminf _{n \rightarrow+\infty} \Big( \int_{\mathbb{R}^{2}} \left(\frac{1}{2}-\frac{1}{\theta}\right)\left(|(\psi_{n})_{x}|^{2} + |D^{-1}_{x} (\psi_{n})_{y}|^{2} \right) dxdy+\int_{\mathbb{R}^{2}}\left(h_{n}^{1}+h_{n}^{2}+h_{n}^{3}\right) d xdy \Big) \\
		& \geq  \int_{\mathbb{R}^{2}} \left(\frac{1}{2}-\frac{1}{\theta}\right)\left(|\psi_{x}|^{2} + |D^{-1}_{x} \psi_{y}|^{2} \right) dxdy+\int_{\mathbb{R}^{2}}\left(h^{1}+h^{3}\right) d xdy=c_{0}.
	\end{aligned}
	$$
	The above inequalities imply that
	\begin{equation}\label{equa110}
		\lim _{n \rightarrow+\infty}  \int_{\mathbb{R}^{2}} \left(|(\psi_{n})_{x}|^{2} + |D^{-1}_{x} (\psi_{n})_{y}|^{2} \right) dxdy=\int_{\mathbb{R}^{2}} \left(|\psi_{x}|^{2} + |D^{-1}_{x} \psi_{y}|^{2} \right) dxdy
	\end{equation}
	and
	$$
	h_{n}^{1} \rightarrow h^{1}, \quad h_{n}^{2} \rightarrow 0 \,\,\text { and }\,\, h_{n}^{3} \rightarrow h^{3} \,\,\,\text {in} \,\, L^{1}\left(\mathbb{R}^{2}\right).
	$$
	Hence,
	$$
	\lim _{n \rightarrow+\infty} \int_{\mathbb{R}^{2}} V(\epsilon_{n} (x_{n}, y_{n}) +\epsilon_{n}(x, y)) |\psi_{n}|^2 d xdy=\int_{\mathbb{R}^{2}} V(x_{0}, y_{0})\psi^{2} dxdy,
	$$
which leads to
	\begin{equation}\label{equa111}
		\lim _{n \rightarrow+\infty} \int_{\mathbb{R}^{2}}\left|\psi_{n}\right|^{2} d xdy=\int_{\mathbb{R}^{2}}|\psi|^{2} dxdy.
	\end{equation}
	As $X$ is a Hilbert space, from \eqref{equa110} and \eqref{equa111},  we can obtain \eqref{strong} and the proof is completed.\;\;\; $\Box$

\section{Regularity} \label{var}
In this section, we study the regularity of the solutions of $(G_{\varepsilon})$ because it is crucial to show that the solutions obtained in the previous section of the modified
problem $(G_{\varepsilon})$ are solutions of the original problem $(P_{\varepsilon})$ as well as to study the concentration phenomenon. The regularity will be obtained by using the Fourier transform of a tempered distribution. Next, we recall the definition of the Fourier transform of a tempered distribution in \cite{Fo}. For a distribution $f$ and multi-index $\alpha$, the derivative of $f$ is given by
$$
\left\langle\partial^{\alpha} f, \phi\right\rangle=(-1)^{|\alpha|}\left\langle f, \partial^{\alpha} \phi\right\rangle,\,\, \text{for any}\,\, \phi \in C_{0}^{\infty}\left(\mathbb{R}^{2}\right) .
$$
The Fourier transform $\widehat{f}$ of a tempered distribution is defined by
$$
\langle\widehat{f}, \phi\rangle=\langle f, \widehat{\phi}\rangle, \,\, \text{for any}\,\, \phi \in \mathcal{S}\, \text { (Schwartz space) },
$$
likewise, the inverse Fourier transform of a tempered distribution $f$ is defined, denoted by $f^{\vee}$, by
$$
f^{\vee}(x)=(2 \pi)^{-2} \widehat{f}(-x).
$$
All basic properties of the usual Fourier transform remain valid for the Fourier transform of a tempered distribution, for instance, for all tempered distribution $f$ and multi-index $\alpha$, we have
$$
\widehat{\partial^{\alpha} f(x)}=i^{|\alpha|} \xi^{\alpha}\widehat{f}(\xi) \,\,\text { and }\,\, \widehat{x^{\alpha} f(x)}=i^{|\alpha|} \partial^{\alpha}\widehat{f}(\xi).
$$
Now, we state and prove the following result.

\begin{theorem}
If $p \in (1,4)$, then any solution $u$ of $(G_{\varepsilon})$ is continuous. Moreover,
	$$
	u \in W^{2, q^{\prime}}\left(\mathbb{R}^{2}\right) \quad \text{with} \left\{
 \begin{array}{l}
 q^{\prime}=\frac{6}{p+1}, \,\,\,\,\,\text { if }\,\, p \neq 3, \\
q^{\prime}\in(1, \frac{3}{2}),\,\, \,\text { if }\,\, p=3.
 \end{array}
 \right.
	$$

	In addition
	$$
	u(x) \rightarrow 0 \,\, \text { as } \,\,|x| \rightarrow \infty .
	$$
\end{theorem}
\noindent\textbf{Proof.} First of all, we must observe that if $u\in X_{\varepsilon}$  is a solution of $(G_\varepsilon)$, then it is a solution, in the distribution sense, of the following problem
\begin{equation}\label{eq9}
	-\Delta u + u_{xxxx} = h_{xx}, \,\,  \mbox{in} \,\,  \R^{2},
\end{equation}
where
$$
h(x,y)=(V(\epsilon x,\epsilon y)-1) u(x,y)-g(\epsilon x, \epsilon y,u).
$$

Since $u \in L^{6}(\R^2)$ and $V-1 \in L^{\frac{6}{p}}(\R^2)$, we have that $(V(\epsilon x,\epsilon y)-1) u(x,y), g(\epsilon x, \epsilon y,u) \in L^{6/(p+1)}(\R^2)$, and so, the  function $h$ also belongs to $L^{6/(p+1)}(\R^2)$ and we can assume that it is a tempered distribution in $\R^2$. Applying the Fourier transform in (\ref{eq9}), in the sense of the tempered distribution,  we have
$$
<\widehat{-\Delta u}, \phi>+<\widehat{u_{xxxx}},\phi>=<\widehat{h _{xx}},\phi>, \,\, \text{for any}\,\, \phi \in \mathcal{S}.
$$

By  the above mentioned properties, we obtain for each $\phi \in \mathcal{S}, $
$$
\begin{array}{c}
	<-\Delta u,\widehat{ \phi}>+<u_{xxxx},\widehat{\phi}>=<h _{xx},\widehat{\phi}>,\\ \\
	< u,-\Delta(\widehat{ \phi})>+  <u,\widehat{\phi}_{xxxx}>= <h ,\widehat{\phi}_{xx}>,\\ \\
	< u,\, \widehat{\vert(x,y)\vert^2\phi}>+ <u,\,\widehat{x^4 \phi}>= -<h,\,\widehat{x^2 \phi}>, \,\, \text{for}\,\,  (x,y) \in \R^2,\\ \\
	< \widehat{u},\, |(\xi_1,\xi_2)|^2 \phi>+ <\widehat{u},\,\xi_{1}^{4}\phi>=-<\widehat{h},\,\xi_{1}^{2}\phi>, \,\, \text{for}\,\, (\xi_1,\xi_2) \in \R^2,
\end{array}
$$
that is
$$
|\xi|^2 \widehat{u}(\xi)+|\xi_1|^4\widehat{u}(\xi)=-|\xi_1|^2\widehat{h}(\xi), \quad \xi=(\xi_1,\xi_2)\in \R^2,
$$
then
$$
\widehat{u}(\xi)=-\Big(\frac{|\xi_1|^2}{|\xi|^2+|\xi_1|^4}\Big)\widehat{h}\equiv - \Phi_1(\xi)\widehat{h}.
$$
Hence
\begin{equation}\label{eq10}u=\overbrace{(-\Phi_1(\xi)\widehat{h})}^{\vee}.\end{equation}

By \cite[Corollary 1]{L}, $\Phi_1$ is a Fourier multiplier on $L^q(\R^2)$ for all $ q \in (1, \infty),$ and so,
\begin{equation}\label{Multiplier}|\overbrace{\Phi_1(\xi)\widehat{f}}^{\vee}|_{q}\leq C |f|_{q}, \,\, \text{for any}\,\, f \in L^q(\R^2),
\end{equation}
where $C$ is a constant which is independent of $f$ and $\Phi_1$.

As $h \in L^{6/(p+1)}(\R^2)$, combining (\ref{eq10}) with   (\ref{Multiplier}), we conclude

\begin{equation}\label{regularu)}
	u \in L^{6/(p+1)}(\R^2).
\end{equation}
Differentiating (\ref{eq9}) with respect to $x$, twice, we have

\begin{equation}\label{eq12}
	-\Delta ( u_{xx})+(u_{xx})_{xxxx}=h_{xxxx} , \quad  (x,y)\in \R^2.
\end{equation}
Applying the Fourier transform in (\ref{eq12}) again,  we get
$$|\xi|^2 \widehat{u_{xx}}(\xi)+|\xi_1|^4\widehat{u}_{xx}(\xi)=|\xi_1|^4\widehat{h}(\xi), \quad \xi=(\xi_1,\xi_2)\in \R^2,
$$
that is,
$$
\widehat{u_{xx}}(\xi)=\Big(\frac{|\xi_1|^4}{|\xi|^2+|\xi_1|^4}\Big)\widehat{h}\equiv \Phi_2(\xi)\widehat{h}
$$
or equivalently
\begin{equation}\label{eq13}
	u_{xx}=\overbrace{(\Phi_2(\xi)\widehat{h})}^{\vee}.
\end{equation}

By \cite[Corollary 1]{L}, $\Phi_2$ is a Fourier multiplier on $L^q(\R^2)$ for $q=6/(p+1)$, then

\begin{equation}\label{regularuxx}u_{xx} \in L^{6/(p+1)}(\R^2).
\end{equation}
Differentiating (\ref{eq9}) with respect to $y$,  we find
\begin{equation}\label{eq16}
	-\Delta ( u_{y}) + (u_{y})_{xxxx} =  h_{xxy} , \quad  (x,y)\in \R^2.
\end{equation}
Similarly, applying the Fourier transform in equation (\ref{eq16}), it follows that
$$
<h_{xxy},\,\widehat{\phi}>=-<h,\,\widehat{\phi}_{xxy}>=-i<\widehat{h},\,\xi_{1}^{2} \xi_{2}\phi>,
$$
we have
$$
|\xi|^2 \widehat{u_{y}}(\xi)+|\xi_1|^4\widehat{u}_{y}(\xi)=-i |\xi_1|^2\xi_2\widehat{h}(\xi), \quad \xi=(\xi_1,\xi_2)\in \R^2,
$$
that is,
$$
\widehat{u_{y}}(\xi)=\Big(\frac{|\xi_1|^2 \xi_2}{|\xi|^2+|\xi_1|^4}\Big)\widehat{(-i)h}\equiv \Phi_3(\xi)\widehat{(-i)h}.
$$
Hence
\begin{equation}\label{eq17} u_{y}=\overbrace{(\Phi_3(\xi)\widehat{(-i)h})}^{\vee}.\end{equation}

By \cite[Corollary 1]{L}, $\Phi_3$ is a Fourier multiplier on $L^q(\R^2)$ for $q=6/(p+1)$, then

\begin{equation}\label{regularuy}u_{y} \in L^{6/(p+1)}(\R^2).
\end{equation}
We claim
\begin{equation}\label{EQ1}
	u_x  \in  L^{12/(p+1)}(\R^2).
\end{equation}

\noindent
Notice $ u \in W^{\overrightarrow{l}}_{\overrightarrow{p}}(\R^2),$ with $\overrightarrow{l}=(2,0)$ and $ \overrightarrow{p}=(6/(p+1), 6/(p+1)). $  Since $u, u_{xx}\in  L^{6/(p+1)}(\R^2),$ we are going to apply  \cite[Theorem 10.2]{BIN} with $\alpha=(1,0),$  $\overrightarrow{q}=(12/(p+1), 12/(p+1)),$  $ l_1=(2,0)$ and $ l_2=(0,0),$
in this case the line containing the points $ (2,0)$ and $(0,0)$ is just $y=0.$
The point $w=\alpha + \frac{1}{\overrightarrow{p}}-\frac{1}{\overrightarrow{q}}=(1,0)+((p+1)/6, (p+1)/6)-((p+1)/12, (p+1)/12)=\Big((p+13)/12, (p+1)/12\Big)$
does not lie on the above line. Then, there exist  positive constants $C_1, C_2>0$ such that
$$
|u_x|_{12/(p+1)}\leq C_1\Big(|u_{xx}|_{6/(p+1)} + |u|_{6/(p+1)}\Big) + C_2 |u|_{6/(p+1)}.
$$
We recall that, up to now,
\begin{equation}\label{Passo1a}
	u, u_{xx}, u_y \in  L^{6/(p+1)}(\R^2) \ \mbox{and} \ u_x  \in  L^{12/(p+1)}(\R^2),  \,\, \text{for any}\,\, p \in (1,4).
\end{equation}
We claim
\begin{equation}\label{Passo1b}
	u \in L^{r}(\R^2)\quad \mbox{for}\quad \left\{
 \begin{array}{l}
  r \in [6/(p+1), \infty], \,\, \mbox{if}\,\, p \neq 3,  \\
r \in (1, \infty),\,\,\,\, \,\,\,\,\,\quad\quad\mbox{if}\,\,  p=3.
 \end{array}
 \right.
\end{equation}
\noindent {\bf Verification for $p\neq 3$:} \,  Notice $ u \in W^{\overrightarrow{l}}_{\overrightarrow{p}}(\R^2),$ with $\overrightarrow{l}=(2,1)$ and $ \overrightarrow{p}=(6/(p+1), 6/(p+1)). $  Since $u_y, u_{xx}\in  L^{6/(p+1)}(\R^2),$ we are going to apply  \cite[Theorem 10.2]{BIN} with $\alpha=(0,0),$  $\overrightarrow{q}=(\infty,\infty),$  $ l_1=(2,0)$ and $ l_2=(0,1).$ In this case the line containing the points $ (2,0)$ and $(0,1)$ is just $y=-\frac{x}{2}+1.$ The point $w=\alpha + \frac{1}{\overrightarrow{p}}-\frac{1}{\overrightarrow{q}}=(0,0)+((p+1)/6, (p+1)/6)-(0, 0)=((p+1)/6, (p+1)/6)$ does not  lie on the above line. By \cite[Theorem 10.2]{BIN}, $ u \in L^{\infty}(\R^2).$ Thus, for $p \neq 3,$ we have $u\in L^{6/(p+1)}(\R^2)\cap L^{\infty}(\R^2),$ from where it follows that
$$
u \in L^r(\R^2),  \,\, \text{for any}\,\, r\in [6/(p+1),\infty] \,\,\,\mbox{if} \,\,\, p \neq 3.
$$

\noindent {\bf Verification for $p= 3$:} \, Notice $ u \in W^{\overrightarrow{l}}_{\overrightarrow{p}}(\R^2),$ with $\overrightarrow{l}=(2,1)$ and $ \overrightarrow{p}=(3/2,3/2). $  Since $u_y, u_{xx}\in  L^{3/2}(\R^2),$ we are going to apply \cite[Theorem 10.2]{BIN} with $\alpha=(0,0),$  $\overrightarrow{q}=(r, r),$  $ l_1=(2,0)$ and $ l_2=(0,1).$
In this case the line containing the points $ (2,0)$ and $(0,1)$ is just $y=-\frac{x}{2}+1.$
The point $w=\alpha + \frac{1}{\overrightarrow{p}}-\frac{1}{\overrightarrow{q}}=(0,0)+(3/2,3/2)-(r, r)=(3/2 -r, 3/2- r)$ does not  lie on the above line, if $ r \in (1,\infty)$. Thereby, by \cite[Theorem 10.2]{BIN}, $ u \in L^{r}(\R^2)$ for all $r > 1.$

In the sequel, setting
$$
\Upsilon(x,y)=h_{xx}=\Big((V(\varepsilon x, \varepsilon y)-1)u-g(\varepsilon x, \varepsilon y,u)\Big)_{xx},
$$
we see that
$$\begin{aligned}
\Upsilon=&\varepsilon^2 V_{xx}(\varepsilon x, \varepsilon y)u+2\varepsilon V_{x}(\varepsilon x, \varepsilon y)u_x+(V(\varepsilon x, \varepsilon y)-1)u_{xx}-\varepsilon^2 \chi_{xx}(f(u)-\tilde{f}(u))-2\varepsilon \chi_x(f'(u)-\tilde{f}'(u))u_x\\
& -\chi f''(u)(u_x)^2-\chi f'(u)u_{xx}-(1-\chi)\tilde{f}''(u)(u_x)^2-(1-\chi)\tilde{f}'(u)u_{xx}.
\end{aligned}$$
We claim
\begin{equation}\label{Pass1F}
	\Upsilon \in L^{q}(\R^2)\quad \mbox{for}\quad
  \left\{
 \begin{array}{l}
q =6/(p+1)\quad  \mbox{if}\,\, p \neq 3, \\
 q \in (1, 3/2)\,\,\quad \,\,\,\mbox{if}\,\, p=3.	
 \end{array}
 \right.
\end{equation}

\noindent {\bf Verification:} \, Since $\chi,\chi_x$ and $\chi_{xx}$ are  bounded, we can drop it. Let us analyze each term of $\Upsilon$. For the terms involving the function $V$, we recall that
$$
u, u_{xx} \in L^{6/(p+1)}(\mathbb{R}^2) \,\, \mbox{and} \,\,  	u_x  \in  L^{12/(p+1)}(\R^2).
$$
Since $V_{xx}, (V-1)\in L^{\frac{6}{p}}(\mathbb{R}^2)$ and $V_x \in L^{\frac{12}{p+1}}(\mathbb{R}^2)$, it follows that
\begin{equation} \label{CV}
V_{xx}u, \,\, V_{x}u_x \,\, \mbox{and} \,\, (V-1)u_{xx} \in  L^{6/(p+1)}(\R^2),
\end{equation}
for $p\neq  3$.

While for $p=3$, a similar argument shows that
$$
V_{xx}u, V_{x}u_x \,\,\mbox{and} \,\, (V-1)u_{xx} \in  L^{q}(\R^2),\,\,\, q \in (1, 3/2).
$$
As
$$
|f(u)-\tilde{f}(u)|\leq C |u|^{p+1}
$$
for some constant $C>0$ and $u \in L^{6}(\R^2)$, it follows that $\vert f(u)-\tilde{f}(u)\vert \in L^{6/(p+1)}(\R^2)$ for $p\neq  3.$  While for $p=3$
$$
\int_{\R^2} |u|^{t (p+1)}dx=\int_{\R^2} |u|^{4t}dx <\infty, \,\, \mbox{for all}\,\, t \geq 1,
$$
in particular,  $ \vert f(u)-\tilde{f}(u)\vert \in L^{q}(\R^2)$ for any $q \in (1,3/2)$.

For the term $f'(u)u_x$, $(f_2)$ leads to
$$
|f'(u) u_x|\leq C |u|^{p}|u_x|.
$$
Note that  $|u|^{p}|u_x|$ belongs to $ L^{6/(p+1)}(\R^2).$ In fact, from (\ref{Passo1b}) since $|u|_{L^{\infty}} \leq C,$ for $p \neq 3,$ then
$$
\int_{\R^2} |u|^{p 6/(p+1)}|u_x|^{6/(p+1)}dx\leq C(\int_{\R^2} |u_x|^{12/(p+1)}dx)^{1/2} <\infty, \ \mbox{if}\  p \neq 3.
$$
For $p=3,$ $ u \in L^{q}(\R^2) $ for all $ q> 1.$ Then
$$
\int_{\R^2} |u|^{3q} |u_x|^q  dx\leq (\int_{\R^2} |u|^{3sq}dx)^{s}(\int_{\R^2} |u_x|^{s'q}dx)^{s'} <\infty
$$
where $s'q=12/(p+1)=3$. Then, $s'=3/q>1 \Leftrightarrow q<3.$ The term $\tilde{f}'(u)u_x$ can be studied in the same way.

For the term $f''(u)( u_x)^2$, note that
$$
|f''(u)( u_x)^2|\leq C |u|^{p-1}|u_x|^2
$$
for some constant $C>0$. We claim that  $|u|^{p-1}|u_x|^2$ belongs to $ L^{6/(p+1)}(\R^2).$ Indeed, from (\ref{Passo1b}) since $u \in L^{\infty}(\R^2)$ for $p \neq 3,$ then
$$
\int_{\R^2} |u|^{6(p-1)/(p+1)}|u_x|^{12/(p+1)}dx\leq C(\int_{\R^2} |u_x|^{12/(p+1)}dx) <\infty.
$$
For $p=3,$ we know that $ u \in L^{q}(\R^2), $ for all $ q> 1.$ Then
$$
\int_{\R^2} |u|^{q 2} |u_x|^{2q}  dx\leq (\int_{\R^2} |u|^{sq2}dx)^{s}(\int_{\R^2} |u_x|^{s'q2}dx)^{s'} <\infty,
$$
where $s'q2=12/(p+1)=3$, that is, $s'=3/2q>1$. Thereby, $s'>1 \Leftrightarrow q<3/2.$ The term $\tilde{f}''(u)( u_x)^2$ can be treated in the same form.

Finally, for the term $f'(u) u_{xx}$, note that
$$
|f'(u) u_{xx}|\leq C |u|^p |u_{xx}|
$$
for some $C>0$. We claim that $ |u|^p \vert u_{xx}\vert\in L^{6/(p+1)}(\R^2).$ In fact, from (\ref{Passo1b}), $u \in L^{\infty}(\R^2)$ for $p \neq 3,$ then
$$
\int_{\R^2} |u|^{p 6/(p+1)}|u_{xx}|^{6/(p+1)}dx\leq C\int_{\R^2} |u_{xx}|^{6/(p+1)}dx<\infty, \ \mbox{if}\ p \neq 3.
$$
For $p=3,$ $ u \in L^{q}(\R^2)$ for all $ q> 1.$ Thus,
$$
\int_{\R^2} |u|^{3q} |u_{xx}|^{q}  dx\leq (\int_{\R^2} |u|^{3sq}dx)^{s}(\int_{\R^2} |u_{xx}|^{s'q}dx)^{s'} <\infty,
$$
where $s'q=6/(p+1)=3/2$, that is, $s'=3/2q>1$. Note that, $s'>1 \Leftrightarrow q<3/2$. The term $\tilde{f}'(u) u_{xx}$ can be treated in the same way.\\

The above analysis proves the Claim.

\vspace{0.5cm}

Differentiating (\ref{eq9}) with respect to $x$ twice, we get

\begin{equation}\label{eq22}
	-\Delta ( u_{xx}) + (u_{xx})_{xxxx} =  \Upsilon_{xx} , \quad  (x,y)\in \R \times \R.
\end{equation}
Applying the Fourier transfom in the equation (\ref{eq22}), as in (\ref{eq10}), we have
$$
\widehat{u_{xx}}(\xi)=-(\frac{|\xi_1|^2}{|\xi|^2+|\xi_1|^4})\widehat{\Upsilon}\equiv -\Phi_1(\xi)\widehat{\Upsilon},
$$
or equivalently,
\begin{equation}\label{eq23}
	u_{xx}=\overbrace{(-\Phi_1(\xi)\widehat{\Upsilon})}^{\vee}.
\end{equation}

By \cite[Corollary 1]{L}, $\Phi_1$ is a Fourier multiplier on $L^q(\R^2)$ for $ q =6/(p+1)$ if  $p \neq  3,$ and $1\leq q < 3/2$ if $p=3.$ Therefore,
\begin{equation}\label{regularuxx1}
	\begin{array}{rl}
		u_{xx} \in L^{6/(p+1)}(\R^2),& \,\, \mbox{if} \,\, p \neq  3,\\
		u_{xx} \in L^{q}(\R^2), \ \forall\,q\in [1, 3/2),&  \,\, \mbox{if} \,\,  p=3.
	\end{array}
\end{equation}
Similarly,
\begin{equation}\label{regularuxx11}
	\begin{array}{rl}
		u_{xxxx}, u_{xxy},u_{yy} \in L^{6/(p+1)}(\R^2),&\,\, \mbox{if} \,\,  p \neq  3,\\
		u_{xxxx}, u_{xxy},u_{yy}  \in L^{q}(\R^2), \ \forall q \in [1, 3/2), & \,\, \mbox{if} \,\, p=3. \\
	\end{array}
\end{equation}

\vspace{0.5cm}

\noindent By (\ref{eq9}),
\begin{equation}\label{eq121}
	-\Delta u + u =  \Theta,  \,\,  \mbox{in} \,\, \R^2
\end{equation}
where
$$
 \Theta=\Upsilon+u-u_{xxxx}.
$$
Since
\begin{equation}\label{Passo1FFF}
 \left\{
 \begin{array}{l}
 \Theta \in L^{6/(p+1)}(\R^2), \quad\quad \,\,\,\,\,\,\,\,\,\,\,\mbox{if} \,\, p\neq  3, \\
 \Theta \in L^{q}(\R^2), \ \forall q \in [1, 3/2), \,\,\mbox{if} \,\, p=3,
 \end{array}
 \right.
\end{equation}
by \cite[Theorem 1]{Ka},
\begin{equation}\label{Passo1FFG}
	u \in W^{2,q}(\R^2) \,\,\text{for}\,\,  \left\{
 \begin{array}{l}
q = 6/(p+1),\,\,\text{if}\,\, p\neq  3, \\
 q\in [1, 3/2),\,\,\,\,\,\, \,\,\text{if}\,\, p=3.	
 \end{array}
 \right.
\end{equation}
Recalling that for any $q>1$, the embedding
$$
W^{2,q}(\R^2)\hookrightarrow C^{0,\alpha}(\overline{\Omega})
$$
is continuous for any smooth bounded domain $\Omega \subset \R^2$ and $0 < \alpha \leq 2 -2/r,$ it follows that $ u \in C(\R^2).$ Moreover, by using the bootstrapping arguments, there exist $0<r_1<r_2$ and $C>0$ such that
$$
\|u\|_{W^{2,q}(B_{r_1}((x, y)))} \leq C(\|\Theta\|_{L^{q}(B_{r_2}((x, y)))}+\|u\|_{L^{q}(B_{r_2}((x, y)))}), \,\, \mbox{for all}\,\, (x, y) \in \R^2.
$$
Using the Sobolev embeddings, there exists $K>0$ independent of $x$ such that
$$
\|u\|_{C(\overline{B}_{r_1}((x, y)))} \leq K(\|\Theta\|_{L^{q}(B_{r_2}((x, y)))}+\|u\|_{L^{q}(B_{r_2}((x, y)))}), \,\, \mbox{for all}\,\, (x, y)\in \R^2.
$$
The last inequality gives
\begin{equation} \label{L1}
	u(x, y) \to 0 \quad \mbox{as} \quad |(x, y)| \to +\infty.
\end{equation}
This completes the proof of the theorem.\\

\section{Proof of Theorem \ref{gio1}}\label{sec4}
In this section we will use the restriction $p \in (1,2)$. Let $\epsilon_{n}\rightarrow 0$ and  $u_{n}\in X_{\epsilon_{n}}$ such that $J_{\epsilon_{n}}(u_{n})=c_{\epsilon_{n}}$
and $J'_{\epsilon_{n}}(u_{n})=0$. From Lemma \ref{lea}, we know that there exists a sequence $(x_{n}, y_{n})\subset \mathbb{R}^2$
such that \linebreak $\psi_{n}(x, y)=u_{n}(x_{n}+x, y_{n}+y)\rightarrow \psi(x, y)$ in $X$.
Moreover, there exists $(x_{0}, y_{0})\in\Lambda$ such that
$$
\lim_{n\rightarrow\infty}\epsilon_{n}(x_{n}, y_{n})=(x_{0}, y_{0})\quad \text{and}\quad V(x_{0}, y_{0})=V_{1}.
$$
Since $\psi_{n}$ satisfies the following equation
$$
-\Delta \psi_{n}+\psi_{n}=g_{n} \,\,\text {in}\,\, \mathbb{R}^{2},
$$
where
$$
g_{n}(x, y)=\psi_{n}-\frac{\partial^{4}}{\partial x^{4}}\psi_{n}-\frac{\partial^{2}}{\partial x^{2}}\psi_{n}+\frac{\partial^{2}}{\partial x^{2}}\Big(V(\epsilon_{n} (x_{n}, y_{n}) +\epsilon_{n}(x, y))\psi_{n}-g\big(\epsilon_{n} (x_{n}, y_{n}) +\epsilon_{n}(x, y), \psi_{n}\big)\Big).
$$
Repeating the same arguments explored in Section 3, we have $\psi_{n} \in W^{2,q}\left(\mathbb{R}^{2}\right)$ for $q=6 /(p+1)$, because $p \not= 3$. Next, the restriction $p \in (1,2)$ is crucial in our approach, because we don't have a good control of the norm of the sequence $((V(\epsilon_{n} (x_{n}, y_{n}) +\epsilon_{n}(x, y))-1))$ in  $L^{\frac{6}{p}}(\mathbb{R}^2)$, which is crucial in the estimate \eqref{CV}. The convergence  $\psi_{n}\rightarrow \psi $ in $X$ combined with the fact that $(V-1)\in L^{\infty}(\R^2)$ and $q \in (2,6)$ ensures that the sequence
$$
h_n(x,y)=(V(\epsilon_n x+\epsilon_n x_n,\epsilon_n y+\epsilon_n x_n)-1) \psi_n(x,y)-g(\epsilon x+\epsilon_n x_n,\epsilon_n y,\psi_n)
$$
converges in $L^{q}(\R^2)$ to
$$
h(x,y)=(V_1-1) \psi(x,y)-f(\psi(x,y)).
$$
Hence, the Fourier multiplier method used in Section 3 permits to prove that the sequences $((\psi_{n})_{x})$, $((\psi_{n})_{xx})$ and $((\psi_{n})_{xxxx})$ are also strongly convergent in $L^{q}(\R^2)$. Now, setting
$$
\Upsilon_n(x,y)=\Big(V(\epsilon_n x+\epsilon_n x_n,\epsilon_n y+\epsilon_n x_n)-1) \psi_n(x,y)-g(\epsilon x+\epsilon_n x_n,\epsilon_n y,\psi_n)\Big)_{xx},
$$
and arguing as in \eqref{Pass1F}, we also derive that
$$
\Upsilon_n \to \Upsilon= \Big(V_1-1) \psi-f(\psi)\Big)_{xx} \,\, \mbox{in} \,\, L^{q}(\R^2).
$$
The convergences above together with the definition of $g$ and the growth of $f$ and $\tilde{f}$ guarantee that
\begin{equation} \label{infinity1}
g_{n} \rightarrow \hat{g} \,\,\,\text{in} \,\,\, L^{q}\left(\mathbb{R}^{2}\right),
\end{equation}
where
$$
\hat{g}(x, y)=\psi(x, y)-\frac{\partial^{4}}{\partial x^{4}}\psi(x, y)-\frac{\partial^{2}}{\partial x^{2}}\psi+V_{1}\frac{\partial^{2}}{\partial x^{2}}\psi(x, y)-\frac{\partial^{2}}{\partial x^{2}}f(\psi(x, y)).
$$
Then, by using again the bootstrapping arguments, there are $0<r_{1}<r_{2}$ and $C>0$ such that
$$
\left\|\psi_{n}\right\|_{W^{2, q}\left(B_{r_{1}}((x, y))\right)} \leq C\Big(\left\|g_{n}\right\|_{L^{q}\left(B_{r_{2}}((x, y))\right)}+\left\|\psi_{n}\right\|_{L^{q}\left(B_{r_{2}}((x, y))\right)}\Big), \,\, \mbox{for all}\,\, (x, y) \in \mathbb{R}^{2} \text { and } n \in \mathbb{N} \text {. }
$$
Using the Sobolev embeddings, there exists $K>0$ independent of $n$ and $(x, y)\in \mathbb{R}^{2}$ such that
$$
\left\|\psi_{n}\right\|_{C\left(\bar{B}_{r_{1}}((x, y))\right)} \leq K(\left\|g_{n}\right\|_{L^{q}\left(B_{r_{2}}((x, y))\right)}+\left\|\psi_{n}\right\|_{L^{q}\left(B_{r_{2}}((x, y))\right)}), \,\, \mbox{for all}\,\, (x, y)\in \mathbb{R}^{2} \text { and } n \in \mathbb{N} \text {. }
$$
The last inequality together with \eqref{infinity1} gives
\begin{equation*} \label{infinity}
\psi_{n}(x, y) \rightarrow 0 \,\, \text { as } \quad|(x, y)| \rightarrow+\infty \,\,\text { uniformly in } n .
\end{equation*}
Now, assume that $\epsilon_{n}(x_{n}, y_{n})=(\hat{x}_{n}, \hat{y}_{n})$. Since $(\hat{x}_{n}, \hat{y}_{n})\rightarrow (x_{1}, y_{1})\in\Lambda$ as $n\rightarrow\infty$, we can find $r>0$ such that $B_{r}((\hat{x}_{n}, \hat{y}_{n}))\subset \Lambda$ for all $n\in \N$.
Therefore $B_{r/\epsilon_{n}}((x_{n}, y_{n}))\subset \Lambda_{\epsilon_{n}}$, $n\in \N$. As a consequence
$$
\mathbb{R}^{2}\backslash \Lambda_{\epsilon_{n}}\subset \mathbb{R}^{2}\backslash B_{r/\epsilon_{n}}((x_{n}, y_{n}))\quad\text{for any}\,\,n\in \N.
$$
By the previous proof, there exists $R>0$ such that
$$
  \vert \psi_{n}(x, y)\vert<a \quad\text{for any}\,\,\vert (x, y)\vert>R, \,\,n\in \N,
$$
where $\psi_{n}(x, y)=u_{\epsilon_{n}}(x+x_{n}, y+y_{n})$. Hence $\vert u_{\epsilon_{n}}\vert<a$ for any $(x,y)\in \mathbb{R}^{2}\backslash B_{R}((x_{n}, y_{n}))$ and $n\in \N$. Then there exists $n_{0}\in \N$ such that for any $n\geq n_{0}$ and $r/\epsilon_{n} >R$ it holds
$$
\mathbb{R}^{2}\backslash \Lambda_{\epsilon_{n}}\subset \mathbb{R}^{2}\backslash B_{r/\epsilon_{n}}((x_{n}, y_{n}))\subset \mathbb{R}^{2}\backslash B_{R}((x_{n}, y_{n})),
$$
which gives $\vert u_{\epsilon_{n}}\vert<a$ for any $(x, y)\in \mathbb{R}^{2}\backslash \Lambda_{\epsilon_{n}}$ and $n\geq n_{0}$.

This means that there exists $\epsilon_{0}>0$,  problem $(P_{\varepsilon})$ has a solution $u_{\epsilon}$ for all $\epsilon\in (0, \epsilon_{0})$. Taking
$v_{\epsilon}(x)=u_{\epsilon}(\frac{1}{\epsilon}(x, y))$, we can infer that $v_{\epsilon}$ is a solution to problem $(P_{\epsilon}^{*})$.

 Finally, we study
the behavior of the maximum points of $\vert v_{\epsilon}(x, y)\vert$. Take $\epsilon_{n}\rightarrow 0$ and $(u_{\epsilon_{n}})$ be a sequence of solutions to problem $(P_\epsilon)$. By the definition of $g$, there exists $\gamma\in (0, a)$ such that
\begin{equation*}
g(\epsilon_{n} x, \epsilon_{n} y,  t)t=f(t)t\leq \frac{V_{0}}{k} t^{2},\quad\text{for all}\,\, (x, y)\in \mathbb{R}^{2},\, 0\leq t\leq\gamma.
\end{equation*}
Using a similar argument above, we can take $R>0$ such that
\begin{equation}\label{5.14}
\Vert u_{\epsilon_{n}}\Vert_{L^{\infty}(B_{R}^{c}((x_{n}, y_{n})))}<\gamma.
\end{equation}
Up to a subsequence, we may also assume that
\begin{equation}\label{5.15}
\Vert u_{\epsilon_{n}}\Vert_{L^{\infty}(B_{R}((x_{n}, y_{n})))}\geq\gamma.
\end{equation}
Indeed, if \eqref{5.15} does not hold, we have that $\Vert u_{\epsilon_{n}}\Vert_{\infty}<\gamma$, and it follows from
$J'_{\epsilon_{n}}(u_{\epsilon_{n}})(u_{\epsilon_{n}})=0$ that
$$
\int_{\mathbb{R}^{2}} \left(|(u_{\epsilon_{n}})_{x}|^{2} + |D^{-1}_{x} (u_{\epsilon_{n}})_{y}|^{2} +V_{0}|u_{\epsilon_{n}}|^{2}\right)dx dy\leq \int_{\mathbb{R}^{2}}g(\epsilon_{n} x, \epsilon_{n} y,  u_{\epsilon_{n}})u_{\epsilon_{n}}dxdy\leq \frac{V_{0}}{k}\int_{\mathbb{R}^{2}} \vert u_{\epsilon_{n}}\vert^{2}dx.
$$
This fact and $k>\frac{\theta}{\theta-2}$ show that $u_{\epsilon_{n}}\equiv 0$ which is a contradiction. Hence \eqref{5.15} holds.

Taking into account \eqref{5.14} and \eqref{5.15}, we can infer that the global maximum points $p_{n}$ of $\vert u_{\epsilon_{n}}\vert$ belongs to $B_{R}((x_{n}, y_{n}))$, that is
$p_{n}=q_{n}+(x_{n}, y_{n})$ for some $q_{n}\in B_{R}((0, 0))$. Recalling that the associated solution of problem $(P_\epsilon)^{*}$ is of the form $v_{n}(x, y)=u_{\epsilon_{n}}((x, y)/\epsilon_{n})$ , we can see that
a  maximum point $\eta_{\epsilon_{n}}$ of $\vert v_{n}\vert$ is $\eta_{\epsilon_{n}}=\epsilon_{n}(x_{n}, y_{n})+\epsilon_{n}q_{n}$. Since $q_{n}\in B_{R}((0, 0))$,
$\epsilon_{n}(x_{n}, y_{n})\rightarrow (x_{1}, y_{1})$ and $V(x_{1}, y_{1})=V_{1}$, from the continuity of $V$, we can conclude that
\begin{equation*}
\lim_{n\rightarrow\infty}V(\eta_{\epsilon_{n}})=V_{1},
\end{equation*}
which concludes the proof of the theorem.
\vspace{0,1cm}

\noindent {\bf Acknowledgments.} C. O. Alves was partially
supported by  307045/2021-8 and Projeto Universal FAPESQ-PB 3031/2021. C. Ji was partially supported by National Natural Science Foundation of China (No. 12171152) and Natural Science Foundation of Shanghai (No. 20ZR1413900).

\vspace{-0,3cm}

\end{document}